\documentclass{article}
\usepackage{authblk}

\usepackage{mathrsfs}  
\usepackage{sidecap}
\usepackage{amsfonts}
\usepackage{amssymb}
\usepackage{amsmath}
\usepackage{geometry}
\newtheorem{theorem}{Theorem}[section]

\newtheorem{problem}[theorem]{Problem}

\usepackage{multirow}
\usepackage{autonum}
\usepackage{graphicx}
\usepackage{tikz}
\usepackage{bm}
\usepackage{siunitx}
\usepackage{gensymb}
\usetikzlibrary{shapes.geometric, arrows}
\usepackage{url}
\tikzstyle{startstop2} = [rectangle, rounded corners, 
minimum width=3cm, 
minimum height=1cm,
text centered, 
draw=black, 
fill=red!30]
\usepackage[pagewise]{lineno}
\tikzstyle{startstop} = [rectangle, rounded corners, 
minimum width=3cm, 
minimum height=1cm,
text centered, 
draw=black, 
fill=orange!30]

\tikzstyle{io} = [trapezium, 
trapezium stretches=true, 
trapezium left angle=70, 
trapezium right angle=110, 
minimum width=3cm, 
minimum height=1cm, text centered, 
draw=black, fill=blue!30]

\tikzstyle{process} = [rectangle, 
minimum width=4cm, 
minimum height=1cm, 
text centered, 
text width=3.8cm, 
draw=black, 
fill=blue!10]
\tikzstyle{process3} = [rectangle, 
minimum width=4cm, 
minimum height=1cm, 
text centered, 
text width=3.8cm, 
draw=black, 
fill=red!10]
\tikzstyle{process4} = [rectangle, 
minimum width=4cm, 
minimum height=1cm, 
text centered, 
text width=3.8cm, 
draw=black, 
fill=red!30]

\tikzstyle{process2} = [rectangle, 
minimum width=4cm, 
minimum height=1cm, 
text centered, 
text width=3.8cm, 
draw=black, 
fill=blue!30]

\tikzstyle{decision} = [diamond, 
minimum width=3cm, 
minimum height=1cm, 
text centered, 
draw=black, 
fill=green!30]
\tikzstyle{arrow} = [thick,->,>=stealth]

\numberwithin{equation}{section}

\newcommand{\R}{{\mathbb R}}

\newcommand{\mR}{{\mathsf R}}
\newcommand{\mT}{{\mathsf T}}

\newcommand{\mI}{{\mathsf I}}
\newcommand{\mA}{{\mathsf A}}

\newcommand{\mL}{{\mathsf L}}

\newcommand{\mD}{{\mathsf D}}

\newcommand{\mSigma}{{\mathsf \Sigma}}

\title{An efficient hierarchical Bayesian method for the Kuopio tomography challenge 2023}

\author[1]{Monica Pragliola\thanks{monica.pragliola@unina.it}}
\author[2]{Daniela Calvetti \thanks{dxc57@case.edu}}
\author[2]{Erkki Somersalo \thanks{ejs49@case.edu}}

\affil[1]{Department of Mathematics and Applications, University of Naples Federico II, Naples, Italy}
\affil[2]{Department of Mathematics, Applied Mathematics, and Statistics, Case Western Reserve University, Cleveland, USA}

\date{}

\begin{document}
\maketitle

\noindent \textbf{Abstract.}
	The aim of Electrical Impedance Tomography (EIT) is to determine the electrical conductivity distribution inside a domain  by applying currents and measuring voltages on its boundary. Mathematically, the EIT reconstruction task can be formulated as a non-linear inverse problem. The Bayesian inverse problems framework has been applied expensively to solutions of the EIT inverse problem, in particular in the cases when the  unknown conductivity is believed to be blocky.  Recently, the Sparsity Promoting Iterative Alternating Sequential (PS-IAS) algorithm, originally proposed for the solution of linear inverse problems, has been adapted for the non linear case of EIT reconstruction \cite{EIT2024} in a computationally efficient manner. Here we introduce a hybrid version of the SP-IAS algorithms for the nonlinear EIT inverse problem, providing a detailed description of the implementation details, with a specific focus on parameters selection. The method is applied to the 2023 Kuopio Tomography Challenge dataset, with a comprehensive report of the running times for the different cases and parameter selections.

\section{Introduction}
Electrical Impedance Tomography (EIT) is a noninvasive imaging technique aimed at estimating the unknown conductivity of the interior of a domain by injecting currents and measuring the voltages at the boundary. The EIT imaging modality, although typically characterized by lower spatial resolution than other tomographic techniques such as computed tomography,  has the advantage of very fast acquisition times, relatively inexpensive hardware and lack of ionizing radiations, making it a versatile tool in many application fields, e.g., in medicine \cite{eit_med1,eit_med2,eit_med3} and in geophysics \cite{eit_geo1,eit_geo2}.

The mathematical formulation of the EIT problem in its forward and inverse paradigm, and the related properties have been widely explored in literature. The aim of the EIT forward problem is to compute the voltages at the boundary corresponding to given   current injection patters when the conductivity map of the inside of the
body is known. The EIT inverse problem, which seeks to recover the unknown conductivity distribution from the currents/voltages, presents all the challenges of an ill-posed problem, amplified by the fact that the forward operator is nonlinear.

The difficulty of reconstructing  meaningful conductivity maps has motivated  a variety of different regularization techniques for the EIT inverse problem. In general, the aim of  regularization approaches is to promote meaningful solutions by penalizing traits that are either unfeasible or not aligned with the expected nature of the sought solution.  In the EIT inverse problem, it is often assumed that the unknown conductivity map is piecewise constant, or, equivalently, that the conductivity gradient distribution is characterized by a sparse structure \cite{borsic2012primal,harhanen2015edge,kaipio2000statistical,gonzalez2016experimental}.

The blocky assumption also applies to the data of the Kuopio Tomography Challenge 2023 (KTC23)\footnote{https//www.fips.fi/KTC2023.php}. The goal of the KTC23 is to provide segmentations of the reconstructed conductivities for a certain number of targets with increasing level of difficulty, each level being represented by a different number of measured boundary voltages available. 

Our approach to the KTC23 is to formulate the EIT reconstruction problem within a hierarchical Bayesian framework with the a priori belief of piecewise constant unknown conductivity. The blocky prior is formulated as a sparsity assumption of  the vector of the increments between  values of the conductivity at adjacent loci in the body domain. Recently \cite{EIT2024}, the authors have generalized a computationally efficient sparsity promoting Bayesian hierarchical model  originally designed for linear inverse problems \cite{calvetti2015hierarchical,calvetti2018bayes,calvetti2019brain,calvetti2020sparsity}, to non-linear scenarios. More specifically, they have extended  the Iterative Alternating Sequential (IAS) algorithm for computing the Maximum a Posteriori (MAP) estimate for the sparsity promoting Bayesian hierarchical model to the EIT inverse problem, and started to analyze its properties. This new algorithm is the basis of our solver for the EIT inverse problems in the context of the KTC23.

The goal of this work is twofold: the first is to advance and solidify the generalization of the Sparsity Promoting IAS (SP-IAS) algorithm to non-linear settings. Following \cite{EIT2024}, we provide an extension of a hybrid version of IAS, originally proposed for linear forward operators in \cite{calvetti2020sparse}, where it was demonstrated to be capable of combining the sparsity enhancement of non-convex optimization with the robust convergence properties of the convex setup. The second goal is to provide a detailed description of how the proposed algorithm can be tailored to the KTC23 dataset. Special attention is given to parameter selection and comparison of computational running times.

The paper is organized as follows. In Section \ref{sec:eit} we discuss the EIT forward and inverse problem, both in continuous and discrete settings. The formulation of the EIT inverse problem together with an analysis of the IAS algorithm is provided in Section \ref{sec:bayes}. In Section \ref{sec:ias} we outline the hybrid IAS algorithm with a detailed analysis of the numerical steps and provide some insight on how to select the parameters. The results of extensive numerical tests with the KTC23 datasets are presented in Section \ref{sec:test}. Finally, we draw some conclusions in Section \ref{sec:concl}.

\section{The EIT forward and inverse model}\label{sec:eit}

Let $\Omega \subset \R^2$ be a bounded simply connected domain with boundary $\partial \Omega$. The electrical conductivity distribution is denoted by function $\sigma\in L^{\infty}(\Omega)$ which is assumed to be strictly positive, i.e., there exists two finite constants $\sigma_1,\sigma_2>0$ such that for each $x\in\Omega$, $\sigma_1\leq\sigma(x)\leq\sigma_2$. In the EIT acquisition setup, injected currents $\{I_{\ell}\}_{\ell=1}^{L}$ are applied through $L$ electrodes $\{e_{\ell}\}_{\ell=1}^{L}$ on the boundary $\partial\Omega$ of the domain, with $e_\ell$ modeled as non-overlapping intervals, i.e. $e_i\cap e_j=\emptyset$, $i,j=1,\ldots,L$, $i\neq j$. The injected currents induce a static electric voltage potential $u:\Omega\to\R$ inside the domain, as well as a voltage vector $U$ at the electrodes.

The behavior of $u$ is governed by the second order elliptic partial differential equation
\begin{equation}
	\label{eq:CEM1}
	\nabla \cdot (\sigma(x)\nabla u(x)) = 0\,,\quad x\in\Omega\,,
\end{equation}
with Neumann-type boundary conditions
\begin{equation}	\label{eq:CEM2}
	\sigma \frac{\partial u}{\partial n}\bigg|_{\Gamma} = 0\,,\quad \int_{e_{\ell}}\sigma \frac{\partial u}{\partial n}dS = I_{\ell}\,,\; \ell=1,\ldots,L\,,
\end{equation}
where $\Gamma = \partial \Omega \setminus \cup_{\ell=1}^L e_{\ell}$. The uniqueness of $u$ can be guaranteed by augmenting the set of conditions in \eqref{eq:CEM2} with the additional requirement
\begin{equation}\label{eq:CEM3}
\left(u+z_{\ell}\sigma\frac{\partial u}{\partial n}\right)\bigg|_{e_{\ell}} = U_{\ell}\,,\quad \ell=1,\ldots,L\,,
\end{equation}
with $U_{\ell}$ representing the voltage at the electrode $e_{\ell}$. The condition above accounts for the formation of a thin, highly resistive layer between the electrode $e_{\ell}$ and the body, characterized by a contact impedance denoted by $z_{\ell}$. 

Equations \eqref{eq:CEM1}-\eqref{eq:CEM3} together with the conditions
\begin{equation}\label{eq:cond_ui}
\sum_{\ell=1}^{L}U_\ell = 0 \,,\qquad \sum_{\ell=1}^{L}I_\ell=0\,,
\end{equation}
constitute what is typically referred to as the Complete Electrode Model (CEM) \cite{somersalo1992existence}. The aim of the EIT forward problem is to determine the potential $u$ in $\Omega$ and $U_{\ell}$ at the electrode $e_{\ell}$, $\ell=1,\ldots,L$, when the currents $\{I_{\ell}\}_{\ell=1}^L$ are applied, and the contact impedance $z$ as well as the conductivity map $\sigma$ are given. Under these conditions, the CEM model is known to have a unique solution that can be found by solving the weak form of the original problem: we refer to \cite{somersalo1992existence} for the details.

The weak formulation allows to naturally introduce a finite element paradigm for the solution of the forward problem in discrete settings. More specifically, consider a triangular tessellation $\mathcal{T}_h=\{K_{\nu}\}_{\nu=1}^{n_t}$  of $\Omega$, with $h>0$ a parameter denoting the mesh size discretization, and let $\{E_\ell\}_{\ell=1}^{L-1}$ be a basis for the space of real $L$-dimensional vectors satisfying the conditions in \eqref{eq:cond_ui}. Then, as shown in details in \cite{EIT2024}, the boundary voltages are given by
\begin{equation}\label{eq:forw}
U = \mR_{\sigma,z}I\,,
\end{equation} 
where $\mR_{\sigma,z}\in\R^{L\times L}$ is the resistance matrix presenting a non-linear dependence on the conductivity $\sigma$ and the contact impedance $z$.

The inverse problem related to the forward formulation in
\eqref{eq:forw} can be stated as follows.

\begin{problem}
Given a frame $\{I^\ell\}_{\ell=1}^{L-1}$ of linearly independent $L$-dimensional currents applied on the boundary of $\Omega$, and the set of corresponding measured voltages $\{U^{\ell}\}_{\ell=1}^{L-1}$,  
\begin{equation}\label{eq:inpb1}
\text{find  }\sigma\text{  such that  }U^\ell = \mR_{\sigma,z}I^\ell+\epsilon^\ell\,,\quad \ell=1,\ldots L-1\,,
\end{equation}
where $\epsilon^\ell$ is a realization of an $L$-variate zero-mean Gaussian random variable. 
\end{problem}

Above, for the sake of definiteness, it is assumed that the set of current patterns form a basis. In practice, the set of  injected current patterns may be an arbitrary set of vectors satisfying the condition (\ref{eq:cond_ui}). Moreover, the voltages are often not measured on every electrode, i.e., the data may consist of a subsample of the full set of voltage data. 
In general settings, one could be also interested in the estimation of the contact impedances which here are assumed to be known.

We begin by assuming that the conductivity can be written as
\begin{equation}
\sigma(x) = \sigma_0 + \delta\sigma(x)\,,
\end{equation}
where $\sigma_0 > 0$ is a presumably known constant, and the perturbation $\delta\sigma$ vanishes at the boundary $\partial\Omega$. In particular, $\sigma\big|_{\partial\Omega} = \sigma_0$.
Given a triangular tessellation $\mathcal{T}_h$, the conductivity $\sigma$ is approximated by a piecewise linear function,
\begin{equation}
\sigma(x) = \sigma_0 + \sum_{\nu=1}^{n}\xi_\nu \varphi_\nu(x)\,,
\end{equation}
where 
$\varphi_\nu$ is a piecewise linear nodal-based Lagrange basis function, i.e., for each node $p_\mu\in\Omega$ of the mesh, 
\begin{equation}
\varphi_\nu(p_\mu) = \delta_{\nu\mu}, \quad \xi_\nu = \delta\sigma(p_\nu),
\end{equation}
and $n$ is the number of interior nodes, as $\delta\sigma$ is assumed to vanish at the boundary $\partial\Omega$. We point out that in \cite{EIT2024}, an element-based discretization of the conductivity was used, while in the current setting, the discretization mesh for the conductivity and the voltage are not coinciding.
Therefore, instead of seeking a function $\sigma\in L^{\infty}(\Omega)$ suitably bounded and satisfying \eqref{eq:inpb1}, the discretized problem is to estimate
the vector $\xi\in\R^{n}$ representing $\sigma$ in the given tessellation.

The Bayesian formulation of the EIT discrete inverse problem that we solve numerically  is obtained by collecting the measurements into a single column vector
\[
b = \left[\begin{array}{c} U^1 \\ \vdots\\ U^{L-1}\end{array}\right] \in \R^m, \quad m = L(L-1),
\] 
and writing the observation model
\begin{equation}\label{eq:inpb2}
b = \left[\begin{array}{c} \mR_{\sigma,z} I^1 \\ \vdots \\  \mR_{\sigma,z} I^{L-1}\end{array}\right] + \left[\begin{array}{c} \varepsilon^1 \\ \vdots\\ \varepsilon^{L-1}\end{array}\right]
= F(\xi,\sigma_0,z) + \varepsilon.
\end{equation}
The computational forward model is parametrized by the  integrals of $\sigma$ over the triangles. Since our main focus is the estimation of the variables $\xi$, to simplify the notation we omit the reference to the dependency on $\sigma_0$ and $z$, and assume $\varepsilon$ to be drawn from a zero mean Gaussian distribution $\mathcal{N}(0,\mSigma)$, with covariance matrix $\mSigma\in\R^{m\times m}$.

\section{Bayesian formulation of the EIT inverse problem}\label{sec:bayes}

In this section, we formulate the nonlinear EIT inverse problem in \eqref{eq:inpb2} in the Bayesian paradigm where all unknown parameters are modeled as random variables.

In the Bayesian setting, the \emph{likelihood} probability density function (pdf) encodes the likelihood of the observed data if the conductivity were known in terms of the noise statistics. Assuming that the noise is additive and follows a zero-mean Gaussian distribution with  covariance matrix $\mSigma$, the likelihood can be written as 
\begin{equation}
\label{eq:lik}
\pi_{b\mid \xi}(b\mid \xi) \propto {\rm exp}\left( -\frac 12 \big\|\mSigma^{-1/2}\big(b - F(\xi)\big)\big\|^2\right).
\end{equation}
where $\|\cdot\|$ denotes the Euclidean norm. The \emph{a priori} information or beliefs on the unknown are described by the \emph{prior} pdf. Motivated by the nature of the Kuopio tomography challenge data, we design a prior accounting for the blocky nature of the target conductivity by expressing a sparsity assumption for the increments in the conductivity $\sigma$, i.e. in the coefficient vector $\xi$.

Let $N$ be the number of edges in the tessellation $\mathcal{T}_h$ connecting nodes such that al least one of them is an interior node. Define a sparse matrix $\mL\in\R^{N\times n}$ with full-column rank mapping the coefficient vector $\xi$ into the vector of increments on the mesh. More precisely, if $E_\ell = \{p_\mu,p_\nu\}$ is an edge connecting the nodes $p_\mu$ and $p_\nu$, we define
\[
 \mL_{\ell,\mu} = \iota(p_\mu), \quad \mL_{\ell,\nu} = -\iota(p_\nu),
\]
where $\iota(p_\kappa) = 1$ for interior nodes $p_\kappa\notin\partial\Omega$ and $\iota(p_\kappa) = 0$ for boundary nodes $p_\kappa\in\partial\Omega$, and other entries in the $\ell$th row vanish. This way, the matrix $\mL$ is sparse, having at most two non-zero entries in each row, and it is easy to verify that $\mL$ is of full rank $n$. We define
\begin{equation}\label{eq:zeta}
\zeta = \mL \xi \in\R^N\,,\quad \xi = \mL^{\dagger}\zeta = (\mL^{\mT}\mL)^{-1}\mL^{\mT}\zeta\,,
\end{equation}
where $\mL^{\dagger}$ is the pseudoinverse matrix of $\mL$.
For a finite element interpretation of the variable $\zeta$, we refer to \cite{bocchinfuso2024adaptive}.


In the Bayesian framework, a prior expressing the sparsity of a random variable can be formulated in the form of a conditionally Gaussian prior of the form
\begin{equation}\label{prior1}
\pi_{\zeta\mid\theta}(\zeta\mid\theta) \propto \frac{1}{(\theta_1\cdots\theta_N)^{1/2}}{\rm exp}\left( - \frac 12 \sum_{j=1}^N \frac{\zeta_j^2}{\theta_j}\right), \quad \zeta \in {\mathcal R}(\mL).
\end{equation}
The conditionally Gaussian prior (\ref{prior1}) is a restriction of the $N$-dimensional Gaussian density ${\mathcal N}(0,\mD_\theta)$ to the range of the matrix $\mL$, where
\[
\mD_\theta = {\rm diag}(\theta_1,\ldots,\theta_N).
\]  
If $\theta_j$ is small, then a priori the component $\zeta_j$ is close to its zero mean, thus the prior model promotes sparsity if only few of the prior variances are of significant size. It has been shown in \cite{calvetti2019hierachical,calvetti2023bayesian} that a computationally convenient suitable choice for sparsity promotion is to  model $\theta$ as a random variable with mutually independent components following a fat-tailed distribution, such as a generalized gamma distribution,
\begin{equation}\label{prior2}
\pi_\theta(\theta) = \prod_{j=1}^N \pi_{\theta_j}(\theta_j), \quad \pi_{\theta_j}(\theta_j) = \frac{|r|}{\Gamma(\beta)\vartheta_j}\left(\frac{\theta_j}{\vartheta_j}\right)^{r\beta -1}
{\rm exp}\left( - \left(\frac{\theta_j}{\vartheta_j}\right)^r\right),
\end{equation}
where $r\neq 0$, $\vartheta_j>0$ is a scale parameter, and $\beta>0$ is a shape parameter, the roles of which will be revisited below. Combining formulas (\ref{prior1}) and (\ref{prior2}), we can write the joint prior of the pair $(\zeta,\theta)$ as
\begin{align}\label{joint prior}\begin{split}
\pi_{\zeta,\theta}(\zeta,\theta) \;{=}\;& \pi_{\zeta\mid\theta}(\zeta\mid\theta)\pi_\theta(\theta) \\
\;{\propto}\;& {\rm exp}\left( - \frac 12 \sum_{j=1}^N \frac{\zeta_j^2}{\theta_j} - \sum_{j=1}^N\left(\frac{\theta_j}{\vartheta_j}\right)^r +\left(r\beta - \frac 32\right)\sum_{j=1}^N \log\frac{\theta_j}{\vartheta_j}\right), \quad \zeta\in {\mathcal R}(\mL).\end{split}
\end{align}
It follows from Bayes' formula that the posterior density can  be written as
\begin{equation}\label{eq:post}
\pi_{\zeta,\theta\mid b}(\zeta,\theta\mid b) = {\rm exp}\big( - G(\zeta,\theta)\big), \quad  \zeta\in {\mathcal R}(\mL),
\end{equation}
where the quantity 
\begin{equation}\label{Gibbs}
G(\zeta,\theta) = \frac 12 \big\|\mSigma^{-1/2}\big(b - F(\mL^\dagger\zeta)\big)\big\|^2  + \frac 12 \sum_{j=1}^N \frac{\zeta_j^2}{\theta_j} + \sum_{j=1}^N\left(\frac{\theta_j}{\vartheta_j}\right)^r  - \left(r\beta - \frac 32\right)\sum_{j=1}^N \log\frac{\theta_j}{\vartheta_j}
\end{equation}
is called the Gibbs energy. 

Often the information encoded in the posterior pdf \eqref{eq:post}-\eqref{Gibbs} is summarized by the maximizer of the distribution, known as the Maximum A Posteriori (MAP) estimator. The main motivation for collapsing the entire posterior into a single vector is to mitigate the computational burden associated to the exploration of the full posterior distribution. In our application the MAP estimate solves
\begin{equation}\label{eq:map}
\{\zeta^*,\theta^*\}\in \underset{\zeta\in{\mathcal R}(\mL),\,\theta\in\R^N}{\rm argmin}G(\zeta,\theta)\,.
\end{equation}
Efficient computational schemes for the computation of the MAP estimate of the posterior (\ref{eq:post}) in the case of linear forward model operators, i.e., $F(\xi)=\mA\xi$,  have been proposed in the literature \cite{calvetti2019hierachical,calvetti2018bayes,calvetti2020sparse}. 
In these works, the authors minimize the corresponding Gibbs functional by means of the Iterative Alternating Sequential (IAS) algorithm, where at each iteration the $\zeta$ and $\theta$ blocks of variables are updated separately by setting
\medskip

\begin{enumerate}
	\item[(1)] $\zeta$-update:
	\begin{equation}
	\zeta^{k+1}\in\arg\min_{\zeta\in\R^N}G(\zeta;\theta^{k})
	\end{equation}
	\item[(2)] $\theta$-update:
	\begin{equation}
	\theta^{k+1}\in\arg\min_{\theta\in\R^N}G(\theta;\zeta^{k+1})
	\end{equation}
	\end{enumerate}
By organizing the computations appropriately, the constraint $\zeta \in {\mathcal R}(\mL)$ can be made automatic.
In the linear case it has been shown that when the hyperprior is a gamma distribution, i.e. $r=1$, the Gibbs energy functional is strictly convex and the IAS iterates  converge to its unique minimizer - see \cite{calvetti2019hierachical}. 

An analytical study of the convexity of the Gibbs functional and corresponding sparsity promotion for generalized gamma hyperpriors \cite{calvetti2020sparse} has inspired the design of hybrid versions of the original IAS. The starting point for hybrid IAS schemes is the observation that when $r\geq 1$, the Gibbs functional is convex in the $\zeta$-domain, while for $0<r<1$, the convexity of $G$ is restricted  to a proper subset of the domain. The price for the greedier sparsity promotion corresponding to $r<1$ is the presence of local minima. To combine the advantages of a unique minimizer and strong sparsity promotion, it was suggested in   \cite{calvetti2020sparsity} to run the IAS with $r=1$ until the iterates approach enough the unique minimizer, then continue the IAS to minimize the Gibbs energy corresponding to a hyperprior with $0<r<1$ with the understanding that the algorithms will stop at a local minimizer. 

Our goal here is to present the implementation details of the recent extension \cite{EIT2024} of the sparsity promoting IAS algorithm for the solution of the nonlinear EIT inverse problems with the (generalized) gamma hyperprior and with the hybrid IAS version, and to systematically test its performance with the 2023 Kuopio tomography challenge dataset.

\section{Sparsity promoting IAS for the nonlinear EIT inverse problem: implementation details}
\label{sec:ias}

We begin this section with the formulation of the hybrid version of the IAS scheme for the nonlinear EIT inverse problem proposed in \cite{EIT2024}, focusing specifically on the role of the parameters. 

As for the case of linear problems, the hybrid IAS for nonlinear problems can be thought of as seeking to solve a sequence of two optimization problems, whose solutions are MAP estimates corresponding to the conditionally Gaussian posteriors with hyperprior distributions characterized by the parameters
\begin{equation}
(r^{(1)},\beta^{(1)},\vartheta^{(1)})\,,\;r^{(1)}=1\,,\text{   and   }(r^{(2)},\beta^{(2)},\vartheta^{(2)})\,,\;0<r^{(2)}< 1\,,
\label{eq:par_hyb}
\end{equation}
with $\beta^{(1)},\beta^{(2)}>0$ and $\vartheta^{(1)},\vartheta^{(2)}\in\R^N$. The properties of the hybrid IAS scheme defined by  the parameters in \eqref{eq:par_hyb} for a linear forward model are well-understood. Since the sequence of iterates for the first objective function converges to the minimizer of a globally convex cost functional, it is reasonable to take the $k$th approximate solution as a reasonable initial guess for the minimization problem with the second objective function, which is only locally convex. In the nonlinear case, the strict convexity of the Gibbs energy functional for the   gamma hyperprior does not hold in general, however, a partial characterization of the convexity properties of the Gibbs functional corresponding to the gamma hyperprior has been established in the case of noiseless measurements: see, e.g.,  \cite[Theorem 7.2]{EIT2024}.

 In spite of the lack of comprehensive theoretical results, the numerical tests presented in Section \ref{sec:test} suggest that the sparsifying properties of the hybrid IAS for the EIT nonlinear problem are similar to those for the linear case. Regardless of the hyperprior, each IAS iteration alternates between the update of the increments $\zeta$ with fixed variances and the update of the variance vector $\theta$ with the increments fixed at the updated values. Since the forward model depends only on $\zeta$, the linearization steps are only performed inside the nonlinear least squares problem to be solved for the update of the increments $\zeta$. 

\begin{figure}[h!]
	\centering
\begin{tikzpicture}[node distance=2cm]

\node (input) [startstop] {1. Input};
\node (param) [process4, below of=input,yshift = -1.5cm] {2. Parameters selection};
\node (param3) [process3, right of =param, xshift = 4cm] {2b. Algorithmic parameters};
\node (param2) [process3, above of=param3] {2a. Subjective parameters};
\node (param4) [process3, below of =param3] {2c. Tailored parameters};
\node (IAS1) [process2, below of = param, yshift = -2cm] {3. Hybrid IAS: phase I};
\node (IAS11) [process, right of = IAS1, xshift = 4cm] {3a. Multiple model linearizations and ${\xi}$-updates};
\node (IAS12) [process, below of = IAS11] {3b. ${\theta}$-update};

\node (IAS2) [process2, below of = IAS1, yshift = -2cm] {4. Hybrid IAS: phase II};
\node (IAS21) [process, right of = IAS2, xshift = 4cm] {4a. Multiple model linearizations and ${\xi}$-updates};
\node (IAS22) [process, below of = IAS21] {4b. ${\theta}$-update};
\node (output) [startstop,below of = IAS2,,yshift = -1cm] {5. Output};

\draw [arrow] (input) -- (param);
\draw [arrow] (param) -- (param2);
\draw [arrow] (param) -- (param3);
\draw [arrow] (param) -- (param4);
\draw [arrow] (param) -- (IAS1);
\draw [arrow] (IAS11.east) -| ++(0.5,0) |- (IAS12.east);
\draw [arrow] (IAS12.west) -| ++(-0.5,0) |- (IAS11.west);
\draw [arrow] (IAS1) -- (IAS11);
\draw [arrow] (IAS21.east) -| ++(0.5,0) |- (IAS22.east);
\draw [arrow] (IAS22.west) -| ++(-0.5,0) |- (IAS21.west);
\draw [arrow] (IAS2) -- (IAS21);
\draw [arrow] (IAS1) -- (IAS2);
\draw [arrow] (IAS2) -- (output);
\end{tikzpicture}
\caption{ Overview of the hybrid IAS algorithm for the EIT nonlinear inverse problem. The apricot boxes contain the input and the output stages. The red and reddish boxes are related to the process of selection of parameters, either in an automatic or manual fashion. The purple and purplish boxes contain the actual body of the algorithm, i.e. the sequence of the two phases and, within each phase, the alternation between the $\zeta$- and the $\theta$-update.}
\label{fig:chart}
\end{figure}
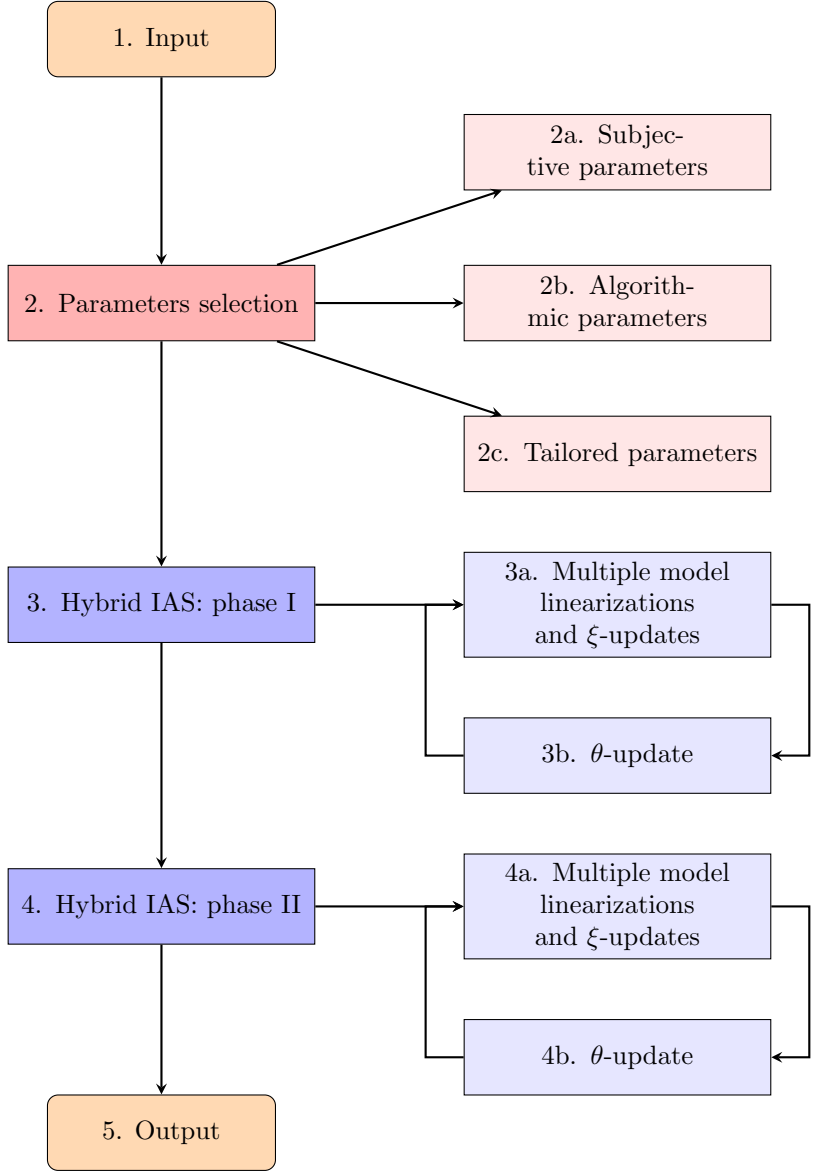

A schematic overview of the algorithm is shown in Figure \ref{fig:chart}. In the following subsections, we discuss in detail how each step can be implemented and how the parameters can be selected for the KTC23 dataset. 

\medskip

\noindent \textbf{1. Input.} In general, the input of hybrid IAS algorithm,  schematically depicted in Figure \ref{fig:chart}, consists of the information related to the acquisition setup, that is the number $L$ of electrodes, their position on the boundary of the domain, the vector of measured voltages $b$, and the vector of contact impedances $z$. In the case of the KTC23 dataset, the conductivity of the water chamber is also available as a reference measurement. Moreover, the mesh used for the discretization of the continuous domain is given as input to our algorithm, together with a  routine that solves the CEM forward model \eqref{eq:CEM1}-\eqref{eq:CEM3}.  In other words, for any given conductivity map and injected current pattern, the solver returns as output the voltages measured at the boundary and the Jacobian of the forward model.

\medskip

\noindent \textbf{2. Parameter selection.} The parameters that need to be assigned can be grouped into three classes, depending on their nature and role. More specifically, the \emph{subjective parameters} are involved in the expression of the hyperprior, the \emph{algorithmic parameters} control the stopping criterion and the number of iterations in the alternating scheme, and, finally, the \emph{tailored parameters}, more specific to the problem to be solved, are mostly related to the definition of the likelihood distribution and to the measurements at hand.

\medskip 

\paragraph{\bf{2a. Subjective parameters}} This class comprises the parameters in \eqref{eq:par_hyb} defining the hybrid hyperprior. In the hybrid IAS paradigm, typically $r^{(1)}=1$, that is the minimization in the first phase of the algorithm is driven by a gamma hyperprior, which, when $\eta^{(1)}=r^{(1)}\beta^{(1)}-3/2\to 0^+$ acts akin to an $\ell_1$-type of penalization. In general, the value of $r^{(2)}$ is chosen in the interval $(0,1)$, although popular choice for this parameter include $r^{(2)}=-1$ - corresponding to an inverse gamma hyperprior for which there is a closed form for the $\theta$-update in step 4b - and $r^{(2)}=1/2$, a choice that that has been proven to work very well in practice in linear scenarios - see, e.g., \cite{calvetti2020sparsity}. Often, instead of selecting a value for $\beta^{(1)}$, it is more natural to assign the value of the related parameter for $\eta^{(1)}$, also referred to as \emph{focality parameter}, according to how sparse the solution is believed to be. In fact, it has been shown that for the linear case the global degree of sparsity of the primary unknown is encoded in the value of the parameter $\eta^{(1)}$: the smaller $\eta^{(1)}$, the sparser the increments vector $\zeta$ is expected to be. The selection of the vector of parameter $\vartheta^{(1)}$, inspired by the role of this vector in the linear case, can be related to the sensitivity of the data to each component of $\zeta$ - see  \cite{calvetti2020sparse,EIT2024}. More specifically, we set
\begin{equation}\label{eq:thstar}
\vartheta_j^{(1)} = \frac{\vartheta^*}{\|\frac{\partial}{\partial \zeta_j} F(\mathsf{L}^{\dagger}\zeta_j)\mid{_{\zeta=0}}\|^2}\,, 
\end{equation}
with $\vartheta^*$ being a positive scalar parameter, and $\frac{\partial}{\partial \zeta_j} F(\mathsf{L}^{\dagger}\zeta_j)\mid_{\zeta=0}$ being the $j$-th column of the Jacobian matrix of the forward model operator computed at $\zeta=0$, that corresponds to the case in which there are no objects immersed in the water chamber. The scaling factor $\vartheta^*$ is manually tuned: in the linear case it may be related to the signal to noise ratio, but we have no information about that in the KTC23 dataset. Finally, the scalar parameter $\beta^{(2)}$ - or, equivalently, $\eta^{(2)}=r^{(2)}\beta^{(2)}-3/2$ - and the vector of parameters $\vartheta^{(2)}$ are automatically determined by the following two conditions introduced to guarantee consistency of the switch between the two hyperpriors. \cite{calvetti2023computationally}:
\begin{enumerate}
	\item Whenever $\zeta_j=0$, the baseline values for $\theta_j$ coincide, which yields
	\begin{equation}\label{eq:conpar1}
	\vartheta^{(1)}\left(\frac{\eta^{(1)}}{r^{(1)}}\right)^{1/r^{(1)}}=	\vartheta^{(2)}\left(\frac{\eta^{(2)}}{r^{(2)}}\right)^{1/r^{(2)}}\,.
	\end{equation}
	\item The marginal expected value for $\theta_j$ is equal using both hypermodels, that is,
	\begin{equation}\label{eq:conpar2}
	\vartheta^{(1)}\frac{\Gamma(\beta^{(1)}+1/r^{(1)})}{\Gamma(\beta^{(1)})} = 	\vartheta^{(2)}\frac{\Gamma(\beta^{(2)}+1/r^{(2)})}{\Gamma(\beta^{(2)})}\,.
	\end{equation}
\end{enumerate}

\medskip

\paragraph{\bf{ 2b. Algorithmic parameters}} The IAS iterations for both hyperpriors can be stopped as soon as either the relative change of the variances at each stage is below a given tolerance, or a maximum number of iterations has been reached. In other words, denoting by $k>0$ the IAS iteration index, we stop iterating as soon as 
\begin{align}\label{eq:rcth}
\begin{split}
\text{I}:&\quad \delta\theta^{k}_{1} = \frac{\|\theta_1^{k+1}-\theta_1^{\text{k}}\|_2}{\|\theta_1^{k}\|_2}<tol\quad\text{or}\quad k\geq k^{(1)}_{\max};\\
\text{II}:&\quad \delta\theta^{k}_{2} = \frac{\|\theta_2^{k+1}-\theta_2^{\text{k}}\|_2}{\|\theta_2^{k}\|_2}<tol\quad\text{or}\quad k\geq k^{(2)}_{\max};
\end{split}
\end{align}
with $\theta_1,\theta_2$ denoting the vector of estimated variances in phase I and II, respectively. Finally, the number of the linearizations and $\zeta$-update in steps 3b, 4b should also be fixed. 

\medskip

\paragraph{\bf{2c. Tailored parameters}} We conclude with a short discussion of the parameters related to the data to be processed, in particular the covariance matrix $\mSigma$ arising in the definition of the likelihood model. Since this information is not part of the KTC23 dataset, we assumed additive zero mean scaled white Gaussian noise, thus letting $\mSigma=\omega^2 \mI$. Due to the lack of information about the noise in the data, $\omega$ has been hand-tuned. Same considerations hold for setting the value of the background conductivity $\sigma_0$.

\medskip

\noindent \textbf{3. Hybrid IAS: phase I} In analogy to the case of linear problems, the first stage of the hybrid IAS is aimed to design a suitable initial guess for the minimization problem addressed in the second stage, which is expected to more strongly promote sparsity in the solution.

\medskip

\paragraph{\bf{3a. Forward model linearization and $\zeta$-update}} Let $k$ denote the iteration index in the first stage of the hybrid IAS. The updated $\zeta^{k+1}$ solves
\begin{equation}\label{eq:subzeta}
\zeta^{k+1} \in\underset{\zeta}{\rm argmin}\left\{\frac{1}{2}\|\mSigma^{-1/2}(b-F(\mL^\dagger\zeta))    \|^2+\frac{1}{2}\|\mD_{\theta^{k}}^{-1/2}\zeta\|^2      \right\}\,,
\end{equation}
with 
\begin{equation}
\mD_{\theta^{k}} = \mathrm{diag}(\theta_1^{k},\ldots,\theta_N^k)\,.
\end{equation}
It follows from our assumption about the statistics of the noise that  $\mSigma^{-1/2}=(1/\omega)\mI$. Introduce the auxiliary variable $\alpha\in\R^N$ defined as
\begin{equation}
\alpha = \mD_{\theta}^{-1/2}\zeta =  \mD_{\theta}^{-1/2}\mL\xi = \mL_{\theta}\xi\,,\quad \text{with}\quad \mL_{\theta}=\mD_{\theta}^{-1/2}\mL\,,
\end{equation}
and 
\begin{equation}
\xi = \mL^\dagger \zeta = \mL_{\theta}^\dagger \alpha\,.
\end{equation}
 Problem \eqref{eq:subzeta} can be reformulated as
\begin{equation}
\alpha^{k+1}\in\underset{\alpha}{\rm argmin}\left\{\|\mSigma^{-1/2}(b-F(\mL_{\theta^k}^\dagger\alpha))\|^2+\|\alpha\|^2\right\}\,.
\label{eq:subalpha}
\end{equation}
Following the derivations in \cite{EIT2024}, we linearize the forward model operator in a neighborhood of $\alpha^{k}=\mD_{\theta^k}^{-1/2}\zeta^k$ to get
\begin{equation}
F(\mL_{\theta^{k}}^\dagger\alpha) \approx  F(\mL_{\theta^{k}}^\dagger\alpha^k) + DF(\mL_{\theta^{k}}^\dagger\alpha^k)\mL_{\theta^{k}}^\dagger (\alpha-\alpha^k) \,,
\end{equation}
where $DF(\mL_{\theta^{k}}^\dagger\alpha^k)$ denotes the Jacobian matrix of the forward model operator with respect to $\xi$ evaluated at $\xi^k=\mL_{\theta^{k}}^\dagger\alpha^k$. 
Defining  
\begin{align}\label{eq:step3a1}
\begin{split}
y \;{=}\;& \mSigma^{-1/2}(b-F(\mL_{\theta^{k}}^\dagger\alpha^k) +DF(\mL_{\theta^{k}}^\dagger\alpha^k)\mL_{\theta^{k}}^\dagger\alpha^k)\\
\mA\;{=}\;& \mSigma^{-1/2} DF(\mL_{\theta^{k}}^\dagger\alpha^k)\mL_{\theta^{k}}^\dagger
\end{split}
\end{align}
we can write \eqref{eq:subalpha} compactly as
\begin{equation}
\alpha^{k+1}\in\underset{\alpha}{\rm argmin}\left\{\|y-\mA\alpha\|^2+\|\alpha\|^2\right\}\,.
\end{equation}
It is straightforward to see that the solution is 
\begin{equation}\label{eq:step3a2}
\alpha^{k+1} = (\mA^\mT\mA+\mI_N)^{-1}(\mA^\mT y)\,,\quad \xi^k = \mL_{\theta^k}^\dagger\alpha^{k+1}\,.
\end{equation}
As pointed out in \cite{EIT2024}, when the dimensionality $m$ of the data is smaller than the number $N$ of increments, as is the case in the present problem, it may be convenient to observe that $\alpha^{k+1} = \mA^\mT w$, where $w$ is the solution of the adjoint problem
\begin{equation}\label{eq:linsys2}
   (\mA \mA^\mT + \mI_m) w = y.
\end{equation}
Summarizing, in Step 3a the update of $\zeta$, i.e. the update of $\xi$, is performed by iterating in an alternating way steps \eqref{eq:step3a1} and \eqref{eq:step3a2} for a certain small number iterations, i.e. 2,3.

\medskip

\paragraph{\bf{3b. $\bm{\theta}$-update}}  The way in which the $\theta$-update is performed does not change when passing from linear to nonlinear settings. Each $\theta_j$ can be updated separately, and the new values must satisfy the first order optimality condition on the cost function, 
\begin{equation}\label{update_theta}
    \theta_j^{k+1}\in\underset{\theta_j}{\rm argmin}\left\{
    \frac 12 \frac{\zeta_j^2}{\theta_j} + \left(\frac{\theta_j}{\vartheta_j}\right)^r  - \left(r\beta - \frac 32\right)\log\frac{\theta_j}{\vartheta_j} 
    \right\}\,.
\end{equation}
Differentiation with respect to $\theta_j$ shows that the updated value of $\theta_j$ solves the one-dimensional non-linear equation
\begin{equation}
    -\frac{1}{2}\left(\frac{\zeta_j}{\theta_j}\right)^2
    (\vartheta_1)_j+r^{(1)}\left(\frac{\theta_j}{(\vartheta_1)_j}\right)^{r^{(1)}-1}-\eta^{(1)}\frac{(\vartheta_1)_j}{\theta_j}=0\,.
\end{equation}
For specific values of $r^{(1)}$, in particular when $r =1$, the solution is available in closed form.

\medskip

\noindent \textbf{4. Hybrid IAS: phase 2} The steps performed at the stage 4a are identical to the ones described above, as the cost functional that is minimized for the $\zeta$-update does not involve the parameters defining the second hyperprior. For what concerns the update of the variances, also in this case each $\theta^{k+1}_j$ can be separately computed by solving a one-dimensional linear equation; a closed form solution is available for $r=r^{(2)}=-1$.

\noindent \textbf{5. Output} In general settings, the output of the hybrid scheme is the vector $\xi$ and the corresponding conductivity map. When applying our approach for the KTC23, an interpolation of the resulting conductivity map is performed on the mesh and a segmentation algorithm is then applied to the final result.

\section{Numerical tests on the KTC2023 dataset}\label{sec:test}

\begin{figure}
	\centering
	\includegraphics[width=6cm]{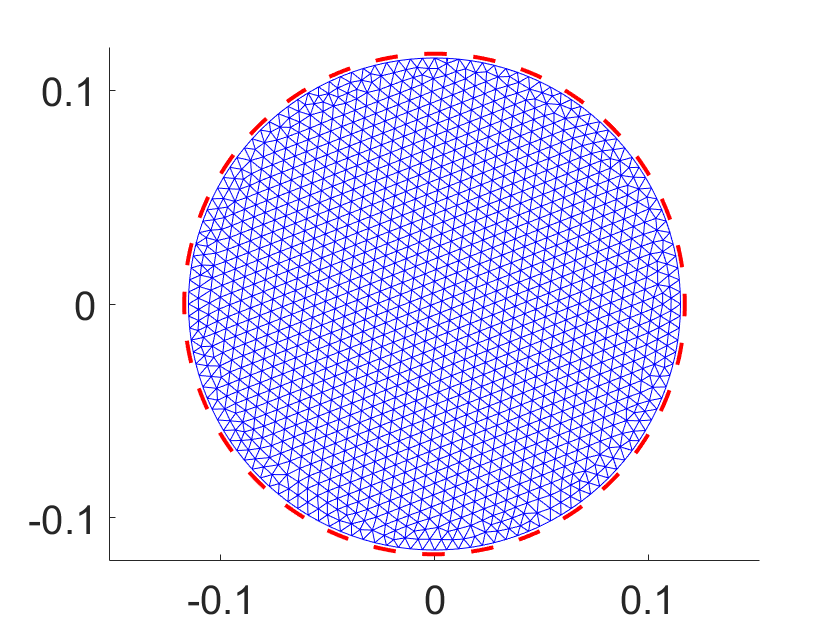}
	\caption{Tessellation used for solving the inverse problem, with the $L=32$ circular arcs representing the areas on the boundary of the circular domain spanned by the electrodes.}
	\label{fig:mesh}
\end{figure}

In this section, we test the sparsity promoting hybrid IAS algorithm for the nonlinear EIT problem discussed in details in Section \ref{sec:ias} on the KTC23 dataset. For the sake of completeness, we describe the experimental set-up used for the data acquisition. Objects with different conductivity/resistivity have been immersed in a circular water chamber, around which $L=32$ electrodes of an EIT device have been installed. The angle spanned by each electrode and the distance between adjacent electrodes is $\alpha=5.625\degree$. For solving the inverse problem, a mesh of $1\,602$ nodes has been provided. In Figure \ref{fig:mesh}, we show the domain tessellation together with the equidistant electrodes.
The number of interior nodes is $1\,473$ which is also the actual dimensionality of the inverse problem, as the values $\xi_j$ corresponding to the boundary nodes are set equal to 0. Also, the entries of the contact impedance vector $z$ are assumed to be all equal to a fixed value $z_0=1\times 10^{-6}$. 

The goal of the KTC23 challenge is to recover a segmented image of the conductivity maps. To this purpose, we first run the hybrid IAS algorithm, then apply the Otsu's method to get the required segmentation. The stopping criterion of the IAS algorithm has been based on the iterations number: more specifically, in \eqref{eq:rcth} we selected $k_{\max}=5$, so that the overall scheme performs a total of 10 iterations. Also, the number of linearizations embedded in the $\zeta$-update is set equal to 2.

The algorithm has been tested on the reconstruction of the conductivity maps induced by 3 different phantoms. For each phantom, we begin by considering the experiment of lowest level of difficulty in which 16 out of 32 electrodes are used for current injections, with the total number of injections $N_{inj}$ being equal to 76, while the voltages between adjacent pairs of electrodes during each current injection is measured. In subsequent tests with increasing level of difficulty, electrodes are progressively removed, so that the number of injected currents decreases: in those cases the quality of the reconstructions is expected to decay. The  number of injections and number of electrodes for each test can be found in Table \ref{tab:1}.

The values of the noise level and of the background conductivity have been hand tuned; in all computed  tests, we set $\omega = 0.004$ and $\sigma_0=0.79$. 

The quality of the reconstruction is assessed by means of a metric based on the structural similarity index (SSIM) \cite{ssim}, taken separately for the conductive and non-conductive inclusions.

All the tests have been performed under Windows 11 in MATLAB R2024a on an Dell PC with an Intel Core i7-13700H @2.4Ghz processor and 16 GB of RAM.

\medskip

\begin{figure}
	\centering
	\begin{tabular}{ccc}
		$r=1$&$r=1/2$&hybrid\\
		\includegraphics[height=3.5cm]{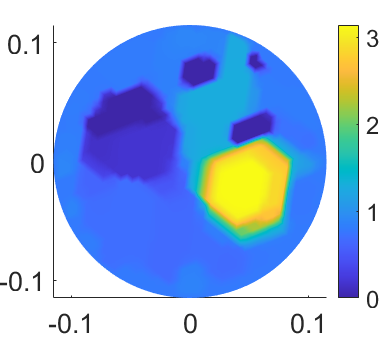}&\includegraphics[height=3.5cm]{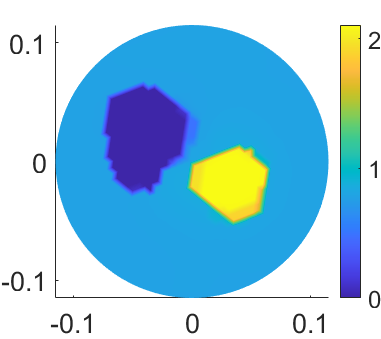}&\includegraphics[height=3.5cm]{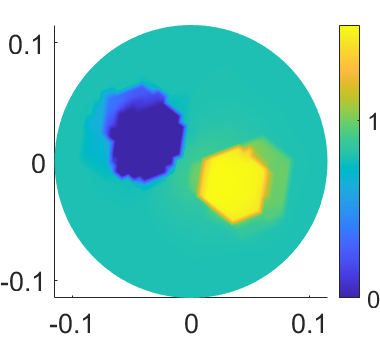}
	\end{tabular}
	\caption{Output reconstruction obtain with the single hyperprior IAS algorithm with $r=1$ (left), $r=1/2$ (middle), and by the hybrid IAS (right) for phantom \#1 and $N_j=76$.}
	\label{fig:hybval}
\end{figure}

\begin{figure}
	\centering
	\setlength{\tabcolsep}{2pt} 
	\begin{tabular}{ccccccc}
		\raisebox{1cm}{\rotatebox{90}{Phase I}}&\multicolumn{2}{c}{\includegraphics[height=3.4cm]{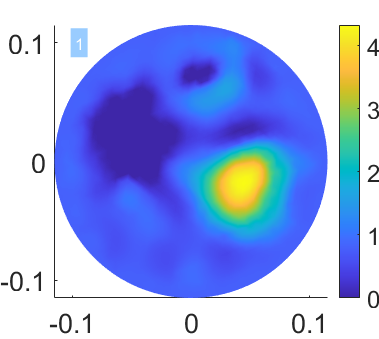}}&\multicolumn{2}{c}{	\includegraphics[height=3.4cm]{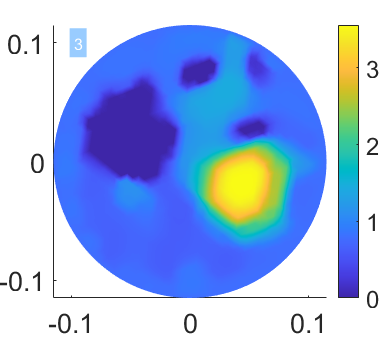}}&\multicolumn{2}{c}{	\includegraphics[height=3.4cm]{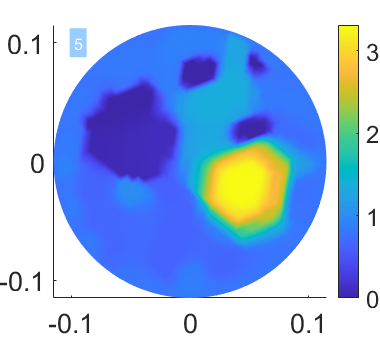}}\\
		\raisebox{1cm}{\rotatebox{90}{Phase II}}&\multicolumn{2}{c}{\includegraphics[height=3.4cm]{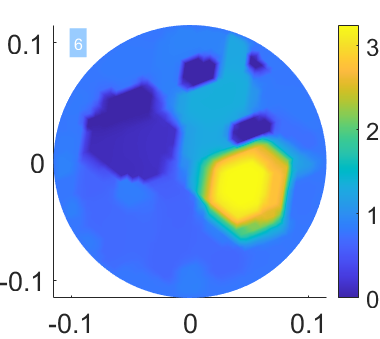}}&\multicolumn{2}{c}{\includegraphics[height=3.4cm]{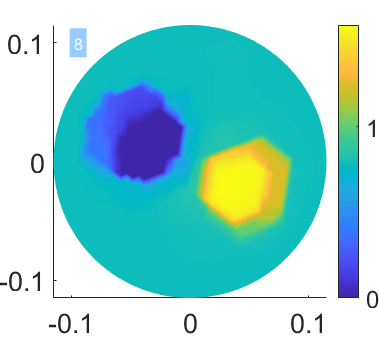}}&\multicolumn{2}{c}{	\includegraphics[height=3.4cm]{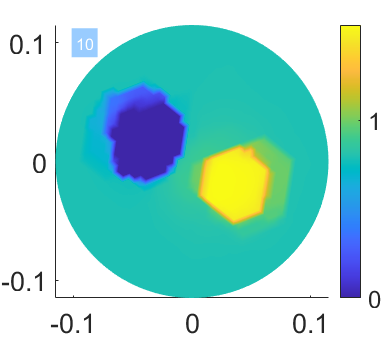}}\\
		&\phantom{xxxxxxxxxxxx}&\multicolumn{2}{c}{\includegraphics[height=2.8cm]{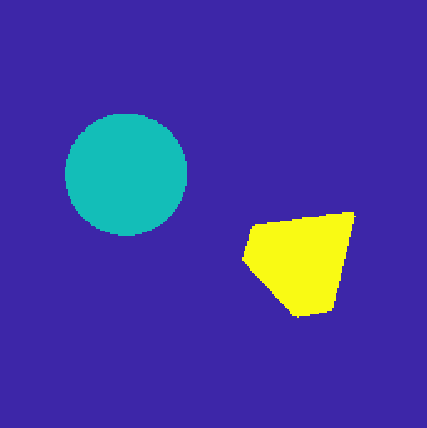}}&\multicolumn{2}{c}{\includegraphics[height=2.8cm]{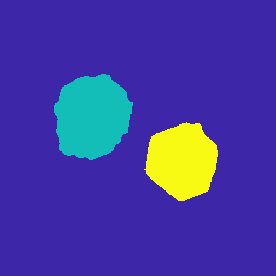}}&
	\end{tabular}
	\caption{Iterative estimates of the conductivity map for phantom \#1 with $N_{inj}=76$ obtained with the hybrid IAS in the first phase with $r^{(1)}=1$ (top panels) and in the second with $r^{(2)}=1/2$ phase (middle panels). In the bottom, target (left) and output segmentation (right) corresponding to the 10-th iteration of the overall hybrid scheme.}\label{fig:out1lev1}
\end{figure}

The first phantom, whose target segmentation is shown in the bottom left panel of Figure \ref{fig:out1lev1}, comprises two inclusions of different shapes. As a preliminary test, we aim to assess the performance of the hybrid IAS scheme with respect to IAS with a single hyperprior. More specifically, for $L=32$ and $N_{inj}=76$, we run the single hyperprior IAS with $r=1$ and $r=1/2$, respectively; we set $\eta^{(1)} = 3 \times 10^{-4}$ and $\vartheta^*$ in \eqref{eq:thstar} equal to $0.03$ , while the parameters identifying the second hyperprior are set according to the conditions \eqref{eq:conpar1}-\eqref{eq:conpar2}. Although our intuition in nonlinear settings is not currently supported by a solid theoretical analysis, we do expect weaker sparsity promotion for $r=1$, and stronger, possibly overbearing sparsification for $r=1/2$. The output reconstructions, after an interpolation on the mesh, are shown in the left and middle panels in Figure \ref{fig:hybval}. One can notice that the gamma hyperprior is not capable of flattening out the background; nonetheless the two objects can be clearly distinguished. On the other hand, the generalized gamma hyperprior with $r=1/2$ produces a solution with sharper edges whose shapes however seem to be less consistent with those of the original targets. Finally, the reconstruction computed by the hybrid IAS scheme, displayed in the right panel, is characterized by a well cleaned background and, in the foreground, by two objects whose dimensions and shapes more accurately resemble the ones of the targets.

As a further analysis, in the top panels of Figure \ref{fig:out1lev1} we show the $\sigma$-estimates during Phase I (first row) and Phase II (second row), with the iteration number being annotated in each subplot. The panels in the second row show how the highly sparsifying characteristics of the functional minimized in Phase II help in flattening the background. The segmentation corresponding to the 10th and last iteration of the hybrid IAS scheme is shown in the bottom right panel. Notice that the polygonal structure of the second highly conductive object seems to be well captured.
%
%
%

\begin{figure}
	\centering
		\setlength{\tabcolsep}{2pt} 
	\begin{tabular}{ccccc}
		\multicolumn{2}{c}{$N_{inj}=56$}&\phantom{x}&	\multicolumn{2}{c}{$N_{inj}=52$}\\
	\includegraphics[height=3cm]{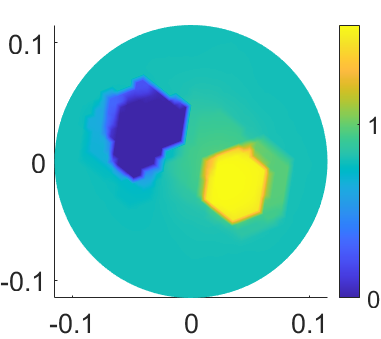}&	\raisebox{0.3cm}{\includegraphics[height=2.5cm]{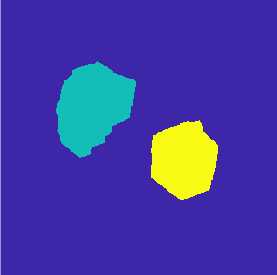}}&\phantom{x}&	\includegraphics[height=3cm]{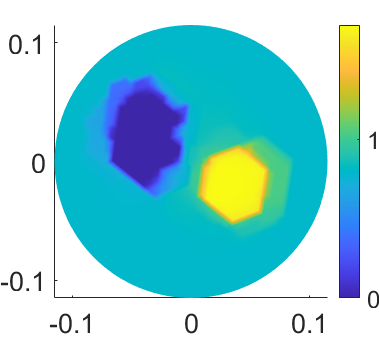}&	\raisebox{0.3cm}{\includegraphics[height=2.5cm]{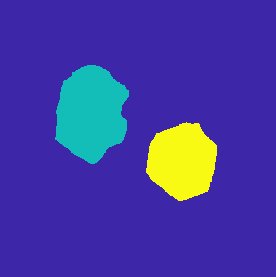}}\\
		\multicolumn{2}{c}{$N_{inj}=48$}&\phantom{x}&	\multicolumn{2}{c}{$N_{inj}=44$}\\
	\includegraphics[height=3cm]{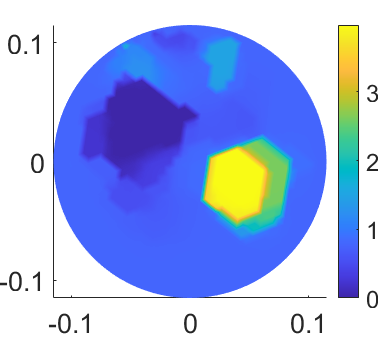}&	\raisebox{0.3cm}{\includegraphics[height=2.5cm]{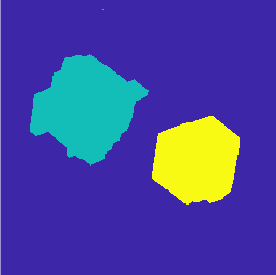}}&\phantom{x}&	\includegraphics[height=3cm]{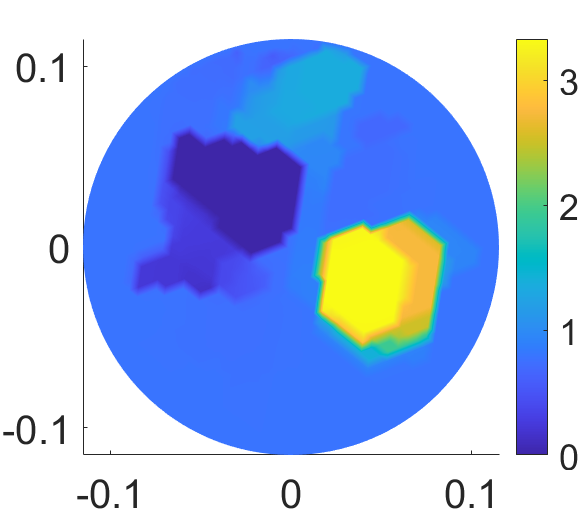}&	\raisebox{0.3cm}{\includegraphics[height=2.5cm]{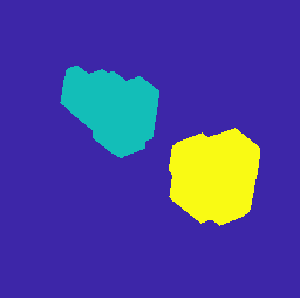}}\\
		\multicolumn{2}{c}{$N_{inj}=30$}&\phantom{x}&	\multicolumn{2}{c}{$N_{inj}=27$}\\
	\includegraphics[height=3cm]{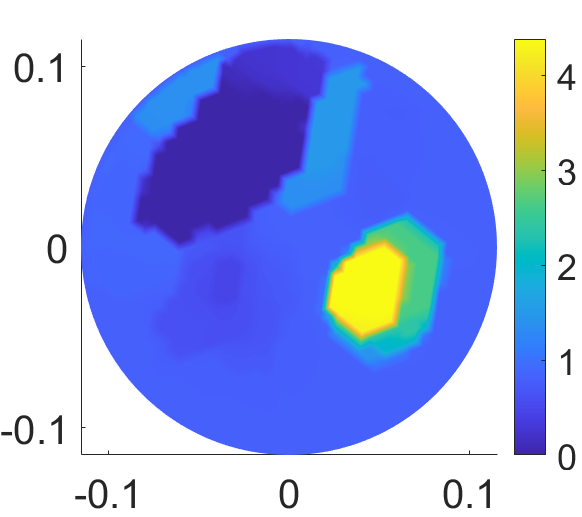}&	\raisebox{0.3cm}{\includegraphics[height=2.5cm]{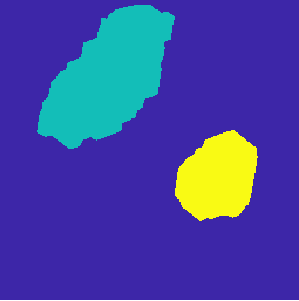}}&\phantom{x}&	\includegraphics[height=3cm]{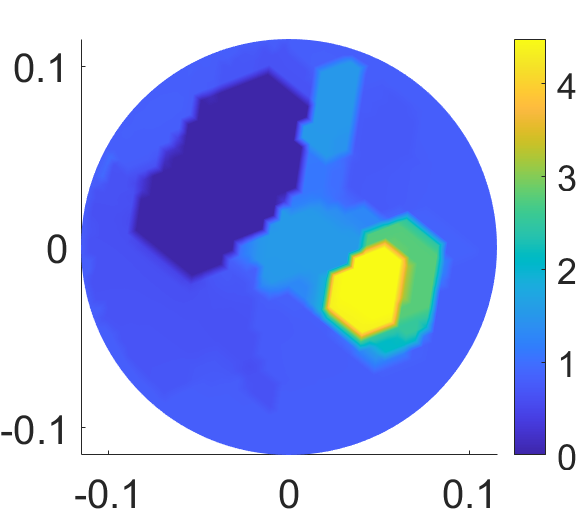}&	\raisebox{0.3cm}{\includegraphics[height=2.5cm]{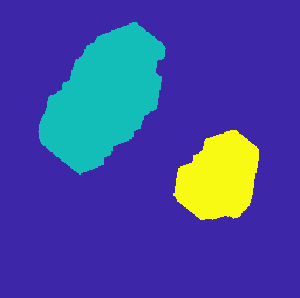}}
	\end{tabular}
\caption{Pairs of output $\sigma$-estimates and corresponding segmentations for phantom \#1 as functions of different numbers $N_{inj}$ of injected currents.}
\label{fig:obj1}
\end{figure}

Figure \ref{fig:obj1} shows the output reconstruction and the corresponding segmentation for $N_{inj}\in\{56,52,48,44,30,27\}$. The value of $\eta^{(1)}$ is selected for the different acquisition scenarios from the interval $[10^{-5},10^{-4}]$, with the focality parameter made smaller as $N_{inj}$ decreases in order to flatten out spurious artifacts generated by the lack of information. The scaling value $\vartheta^*$ has been chosen from the interval $[0.03,0.5]$, with the larger values corresponding to the case of fewer injections. Comparison of the segmented reconstruction with the target in Figure \ref{fig:out1lev1} shows that for $N_{inj}\geq 44$ the dimensions of the phantoms and the profile of the polygonal object are preserved, while for $N_{inj}=30,27$ - corresponding to the case when  22 and 20 electrodes, respectively, are employed - the shapes of the objects seem to be less sharp, although still clearly distinguishable.
The values of the achieved SSIM scores are reported in Table \ref{tab:1}.
\begin{table}
	\centering
	\begin{tabular}{c|c|c|c|c}
	$L$&$N_{inj}$&phantom \#1&phantom \#2&phantom \#3\\
			\hline
32&	76&0.6915&0.8978&0.7628\\
30&	56&0.7031&0.8987&0.7908\\
28&	52&0.6981&0.8939&0.7912\\
26&	48&0.6308&0.8774&0.7651\\
24&	44&0.5582&0.8987&0.8093\\
22&	30&0.5781&0.6978&0.7206\\
20&	27&0.6361&0.6341&0.6317
	\end{tabular}
\medskip
\caption{Values of the SSIM-based scores achieved for the three different phantoms and for different numbers $N_{inj}$ of injective currents, i.e. different number $L$ of electrodes.}
\label{tab:1}
\end{table}
%
%
%

\begin{figure}
	\centering
	\setlength{\tabcolsep}{2pt} 
	\begin{tabular}{ccccccc}
		\raisebox{1cm}{\rotatebox{90}{Phase I}}&\multicolumn{2}{c}{\includegraphics[height=3cm]{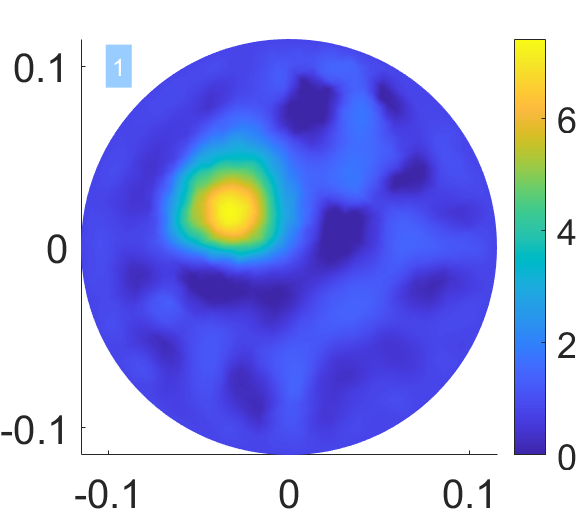}}&\multicolumn{2}{c}{	\includegraphics[height=3cm]{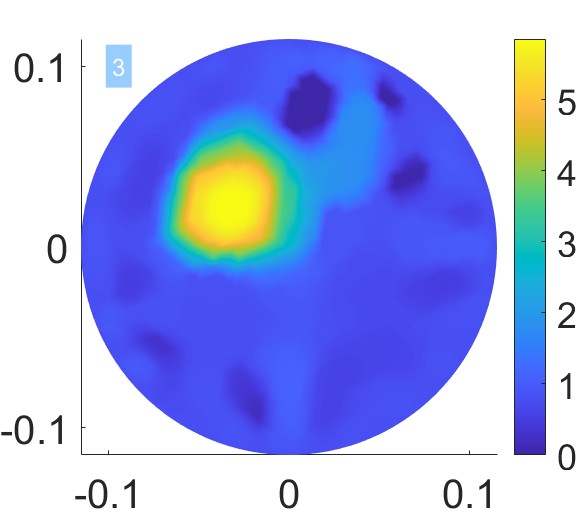}}&\multicolumn{2}{c}{	\includegraphics[height=3cm]{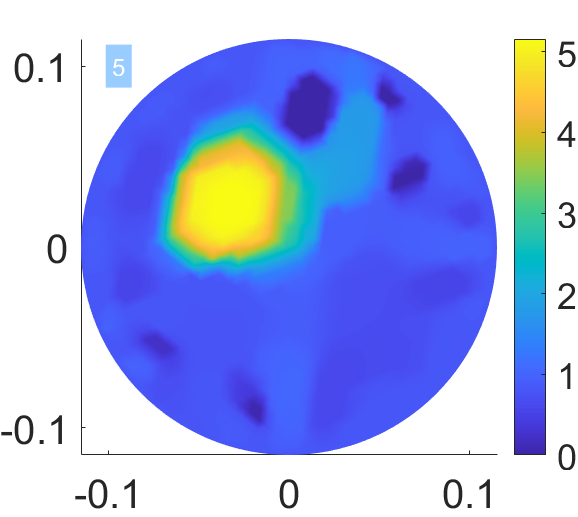}}\\
		\raisebox{1cm}{\rotatebox{90}{Phase II}}&\multicolumn{2}{c}{\includegraphics[height=3cm]{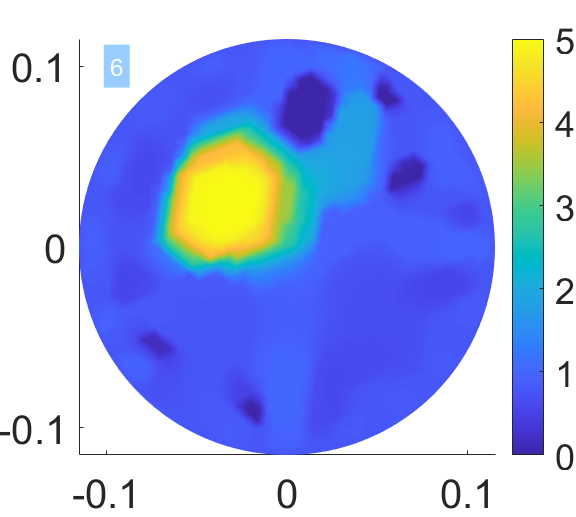}}&\multicolumn{2}{c}{\includegraphics[height=3cm]{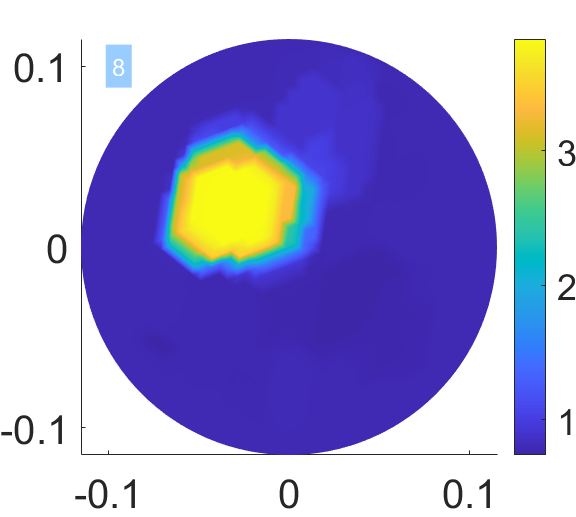}}&\multicolumn{2}{c}{	\includegraphics[height=3cm]{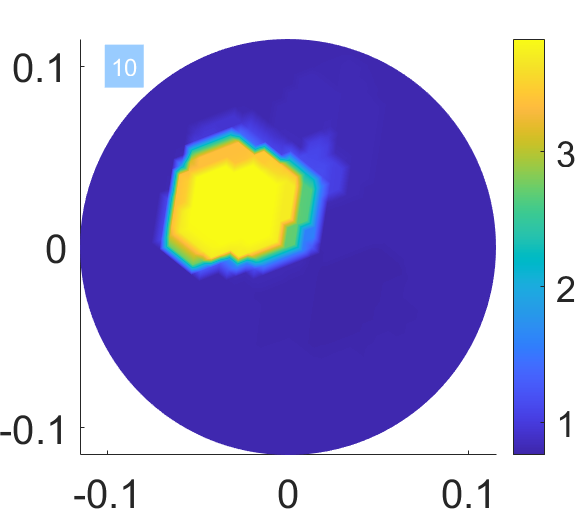}}\\
		&\phantom{xxxxxxxxxxxx}&\multicolumn{2}{c}{\includegraphics[height=2.8cm]{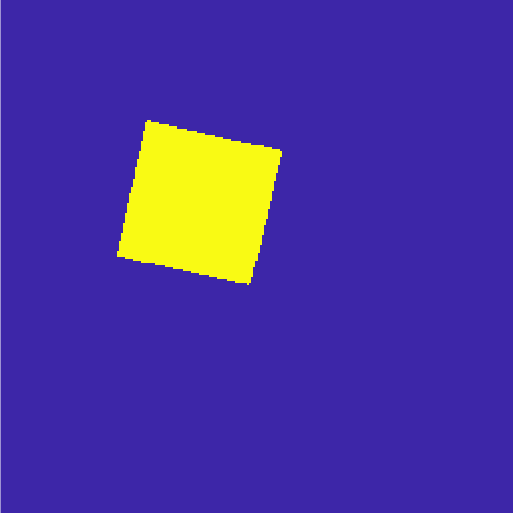}}&\multicolumn{2}{c}{\includegraphics[height=2.8cm]{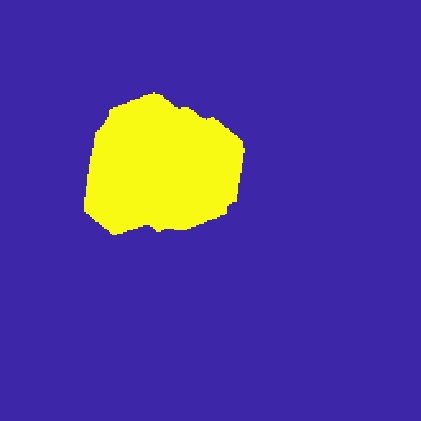}}&
	\end{tabular}
	\caption{Iterative estimates of the conductivity map for phantom \#2 with $N_{inj}=76$ computed with the hybrid IAS scheme in the first phase with $r^{(1)}=1$ (top panels) and in the second with $r^{(2)}=1/2$ phase (middle panels). The bottom row displays the target (left) and output segmentation (right) corresponding to the 10-th iteration of the hybrid scheme.}\label{fig:out2lev1}
\end{figure}

Next we consider the reconstruction of the conductivity map induced by the second phantom, whose target segmentation is shown in the bottom left panel of Figure \ref{fig:out2lev1}. In the top rows of the same figure we show the behavior of the two phases of the hybrid IAS scheme with parameters $\eta^{(1)}=5\times 10^{-6}$ and $\vartheta^*=0.4$ for $N_{inj}=76$: as in the previous case, the  last iteration is characterized by a drastic flattening out of the background. From the output segmentation shown in the bottom right panel, one can conclude that the algorithm is capable of preserving a partial profile and the orientation of the object.

\begin{figure}
	\centering
		\setlength{\tabcolsep}{2pt} 
	\begin{tabular}{ccccc}
		\multicolumn{2}{c}{$N_{inj}=56$}&\phantom{x}&	\multicolumn{2}{c}{$N_{inj}=52$}\\
		\includegraphics[height=3cm]{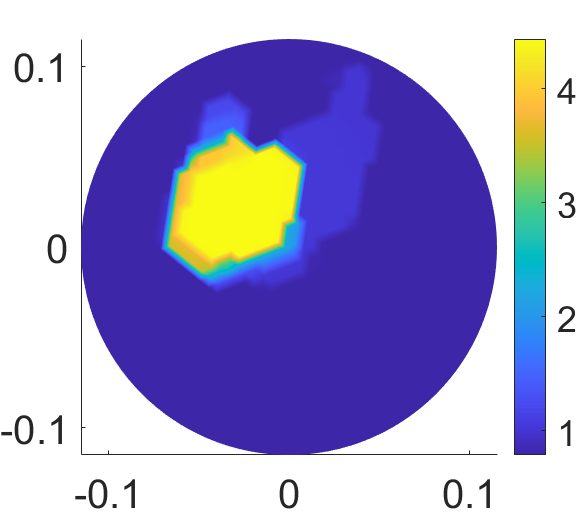}&	\raisebox{0.3cm}{\includegraphics[height=2.5cm]{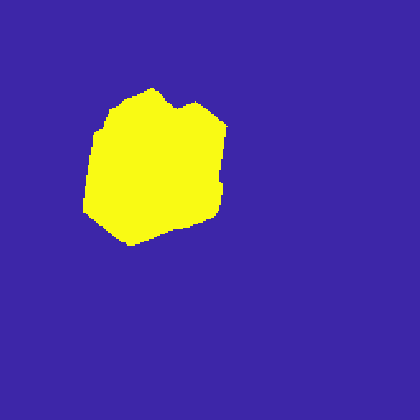}}&\phantom{x}&	\includegraphics[height=3cm]{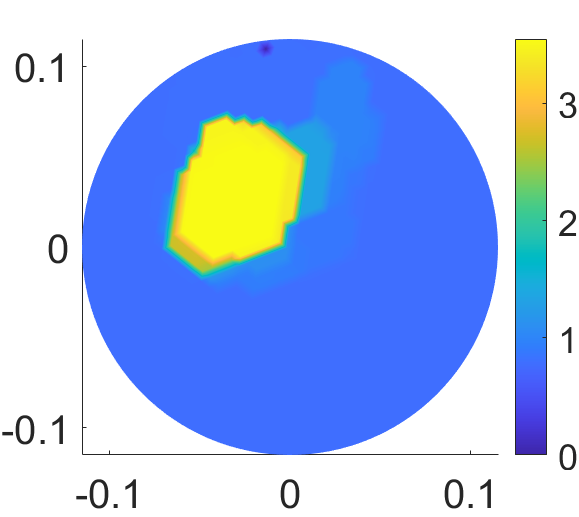}&	\raisebox{0.3cm}{\includegraphics[height=2.5cm]{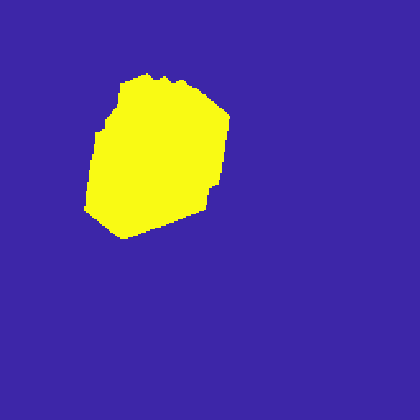}}\\
		\multicolumn{2}{c}{$N_{inj}=48$}&\phantom{x}&	\multicolumn{2}{c}{$N_{inj}=44$}\\
		\includegraphics[height=3cm]{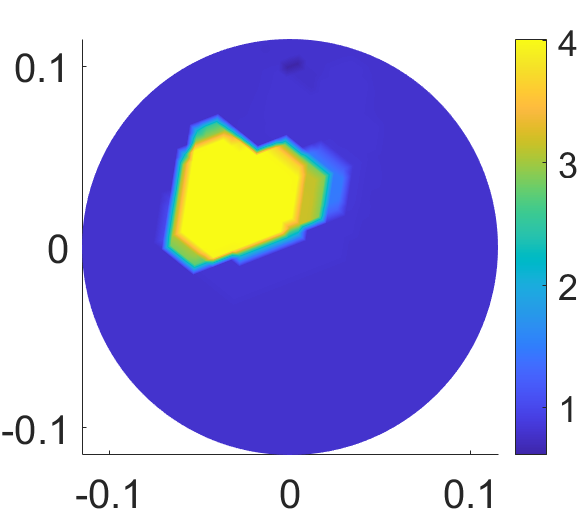}&	\raisebox{0.3cm}{\includegraphics[height=2.5cm]{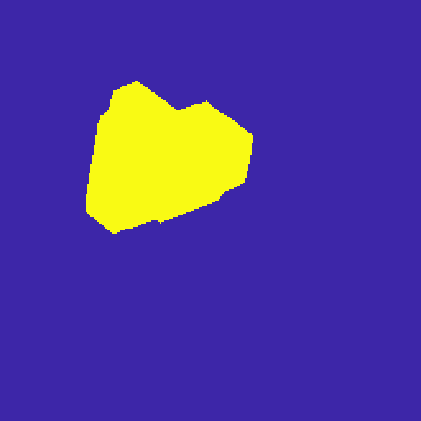}}&\phantom{x}&	\includegraphics[height=3cm]{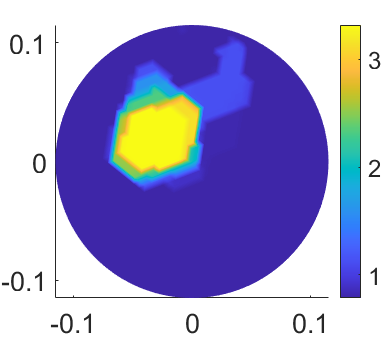}&	\raisebox{0.3cm}{\includegraphics[height=2.5cm]{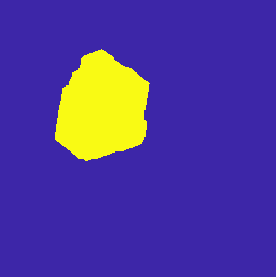}}\\
		\multicolumn{2}{c}{$N_{inj}=30$}&\phantom{x}&	\multicolumn{2}{c}{$N_{inj}=27$}\\
		\includegraphics[height=3cm]{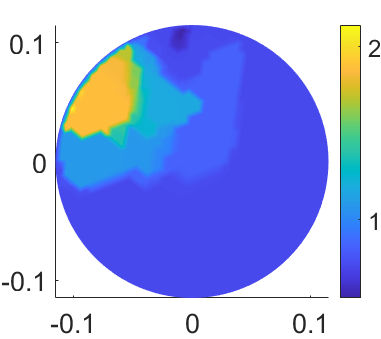}&	\raisebox{0.3cm}{\includegraphics[height=2.5cm]{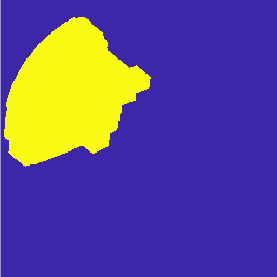}}&\phantom{x}&	\includegraphics[height=3cm]{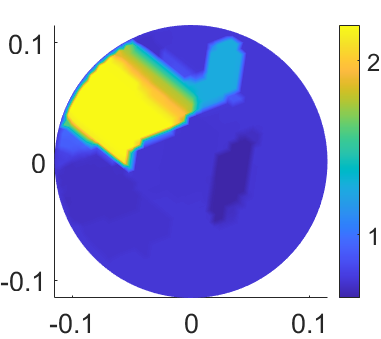}&	\raisebox{0.3cm}{\includegraphics[height=2.5cm]{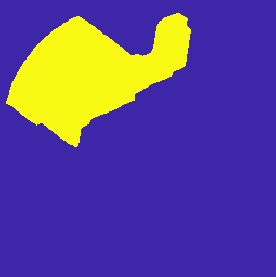}}
	\end{tabular}
\caption{Pairs of output $\sigma$-estimates and corresponding segmentations for phantom \#2 as functions of different numbers $N_{inj}$ of injected currents.}
\label{fig:obj2}
\end{figure}

The output reconstruction and the related segmentation obtained for different values of $N_{inj}$, with parameters $\eta^{(1)},\vartheta^*$ ranging in the intervasls $[10^{-7},10^{-6}]$ and $[0.4,0.6]$, respectively, are shown in Figure \ref{fig:obj2}, while the scores achieved are reported in the fourth column of Table \ref{tab:1}. 

\begin{figure}
	\centering
	\setlength{\tabcolsep}{2pt} 
	\begin{tabular}{ccccccc}
		\raisebox{1cm}{\rotatebox{90}{Phase I}}&\multicolumn{2}{c}{\includegraphics[height=3cm]{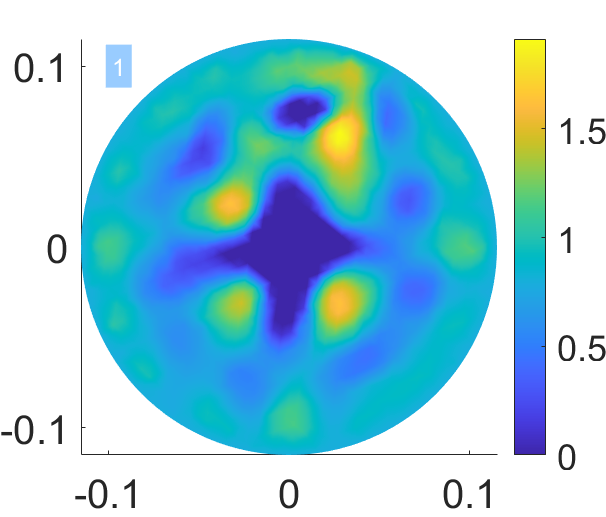}}&\multicolumn{2}{c}{	\includegraphics[height=3cm]{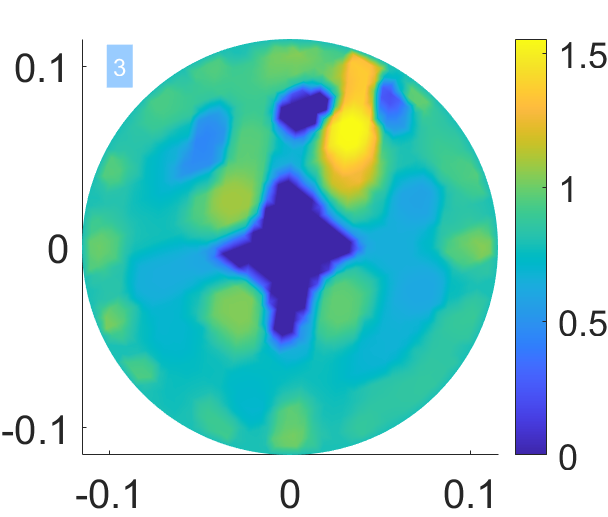}}&\multicolumn{2}{c}{	\includegraphics[height=3cm]{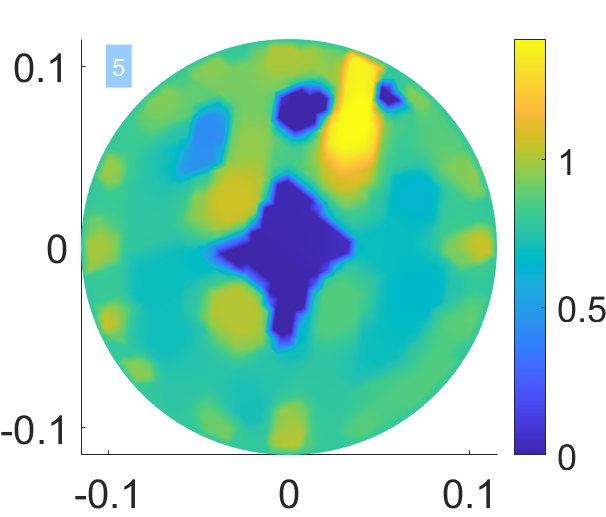}}\\
		\raisebox{1cm}{\rotatebox{90}{Phase I}}&\multicolumn{2}{c}{\includegraphics[height=3cm]{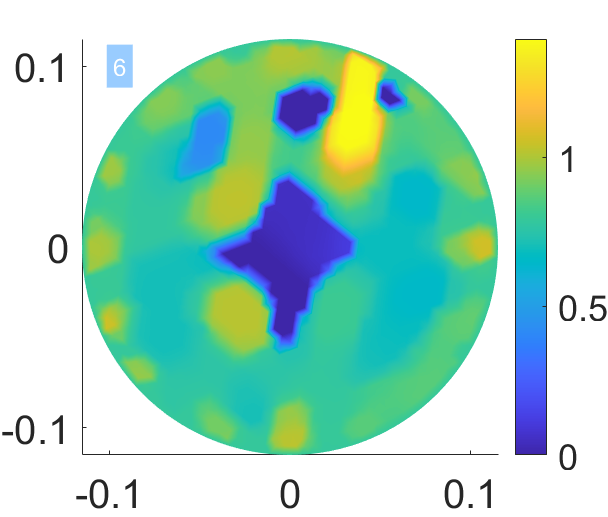}}&\multicolumn{2}{c}{\includegraphics[height=3cm]{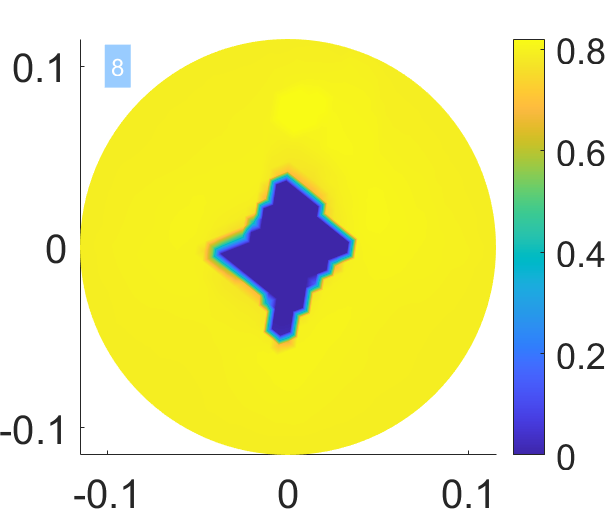}}&\multicolumn{2}{c}{	\includegraphics[height=3cm]{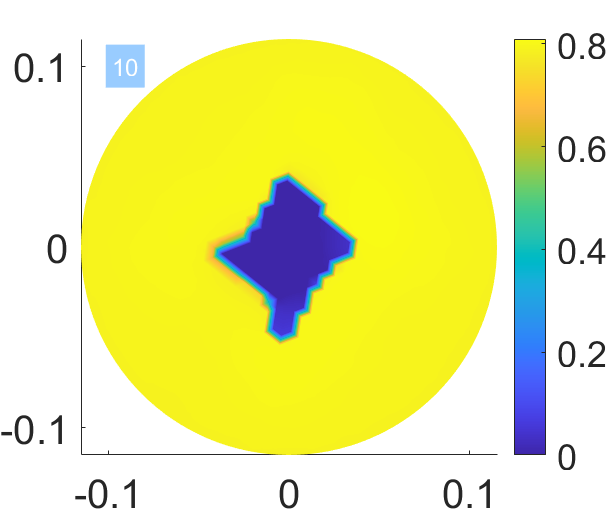}}\\
		&\phantom{xxxxxxxxxxxx}&\multicolumn{2}{c}{\includegraphics[height=2.8cm]{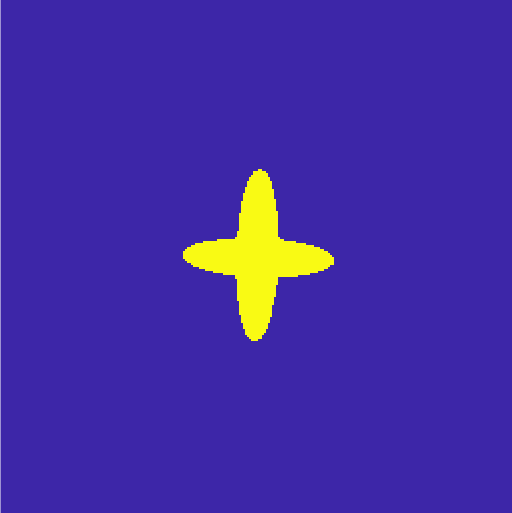}}&\multicolumn{2}{c}{\includegraphics[height=2.8cm]{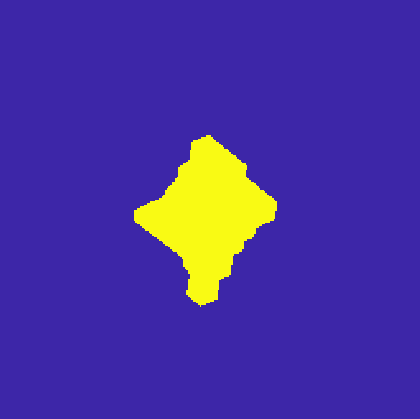}}&
	\end{tabular}
	\caption{Iterative estimates of the conductivity map for phantom \#3 with $N_{inj}=76$ computed by  the hybrid IAS scheme in the first phase with $r^{(1)}=1$ (top panels) and in the second with $r^{(2)}=1/2$ phase (middle panels). In the bottom, target (left) and output segmentation (right) corresponding to the 10-th iteration of the overall hybrid scheme.}\label{fig:out3lev1}
\end{figure}

We finally consider the third phantom, with target segmentation shown in the bottom left panel of Figure \ref{fig:out3lev1}. In this case the inclusion is near the center of the water chamber, which may mitigate the artifacts that tend to arise when the object is close to the boundary. The concave shape of the inclusion, however, poses a challenge, especially when considering the uniform meshes that we use for the reconstructions. For $N_{inj}=76$, the hybrid IAS was run with parameters $\eta^{(1)}=5\times 10^{-4}$ and $\vartheta^* = 0.05$. The iteratively estimated conductivity maps are shown in the first two rows of Figure \ref{fig:out3lev1}. The final reconstruction and the corresponding segmentation confirm that the background is correctly flattened out, while the shape of the object is very close to the target.

\begin{figure}
	\centering
			\setlength{\tabcolsep}{2pt} 
	\begin{tabular}{ccccc}
		\multicolumn{2}{c}{$N_{inj}=56$}&\phantom{x}&	\multicolumn{2}{c}{$N_{inj}=52$}\\
		\includegraphics[height=3cm]{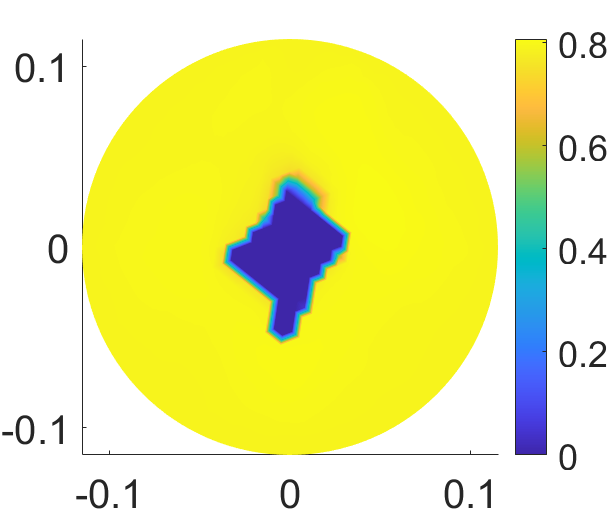}&	\raisebox{0.3cm}{\includegraphics[height=2.5cm]{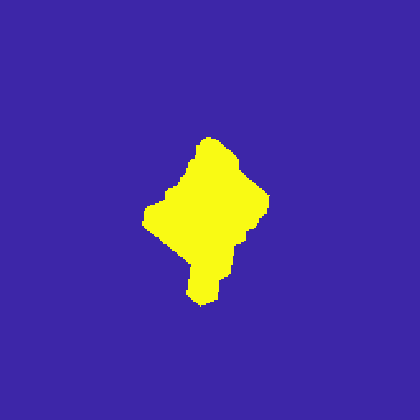}}&\phantom{x}&	\includegraphics[height=3cm]{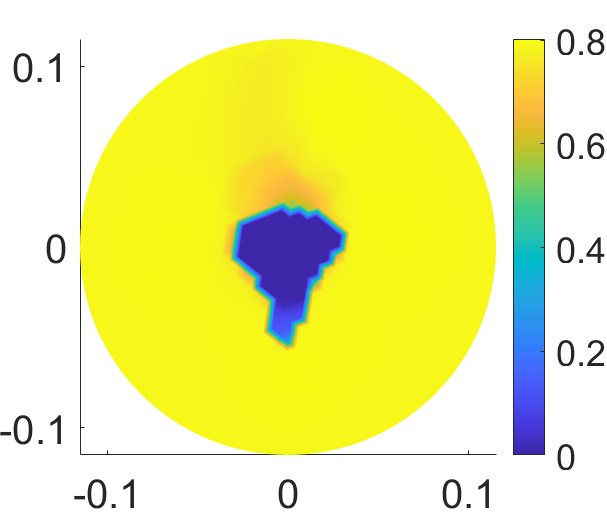}&	\raisebox{0.3cm}{\includegraphics[height=2.5cm]{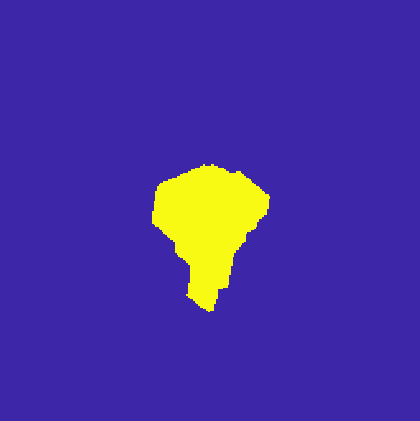}}\\
		\multicolumn{2}{c}{$N_{inj}=48$}&\phantom{x}&	\multicolumn{2}{c}{$N_{inj}=44$}\\
		\includegraphics[height=3cm]{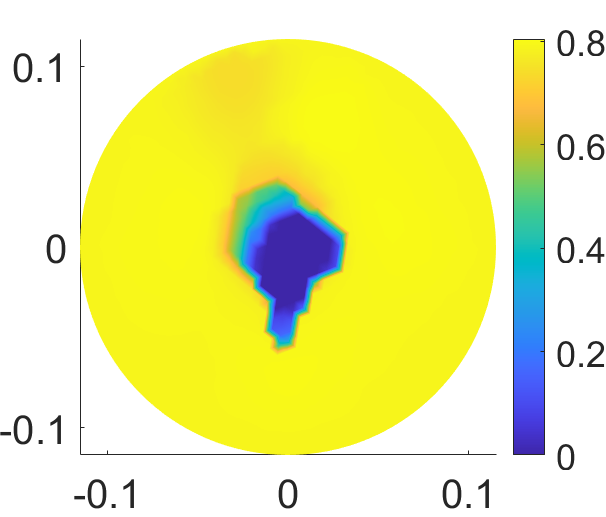}&	\raisebox{0.3cm}{\includegraphics[height=2.5cm]{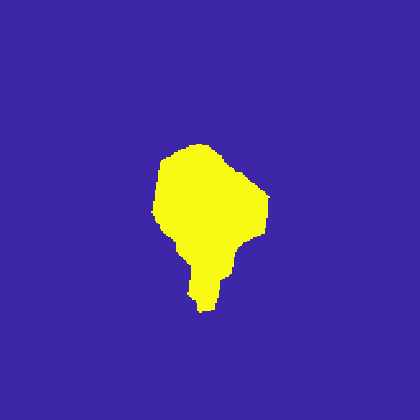}}&\phantom{x}&	\includegraphics[height=3cm]{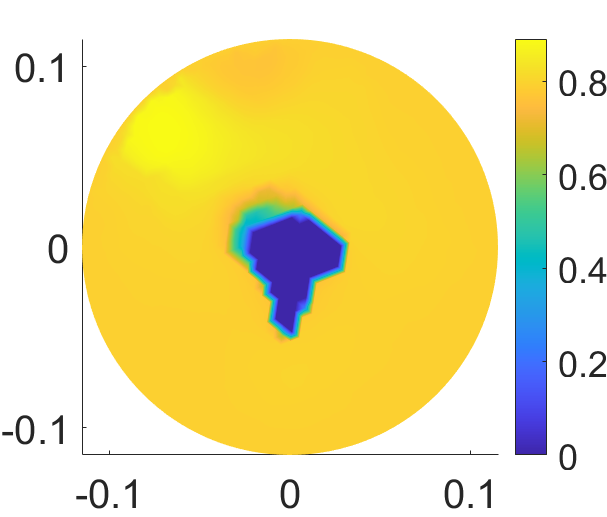}&	\raisebox{0.3cm}{\includegraphics[height=2.5cm]{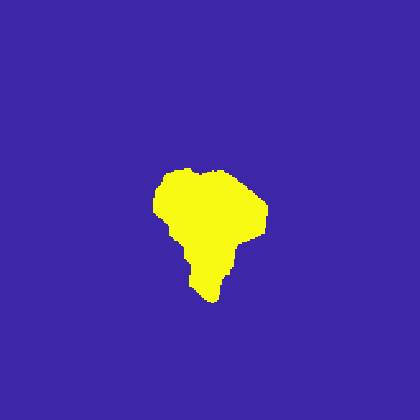}}\\
		\multicolumn{2}{c}{$N_{inj}=30$}&\phantom{x}&	\multicolumn{2}{c}{$N_{inj}=27$}\\
		\includegraphics[height=3cm]{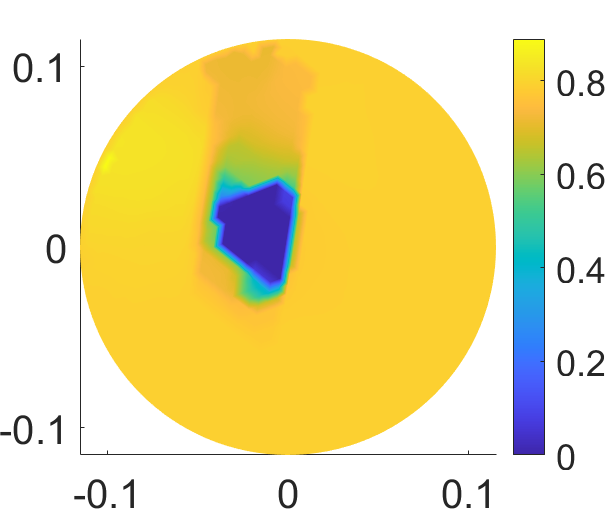}&	\raisebox{0.3cm}{\includegraphics[height=2.5cm]{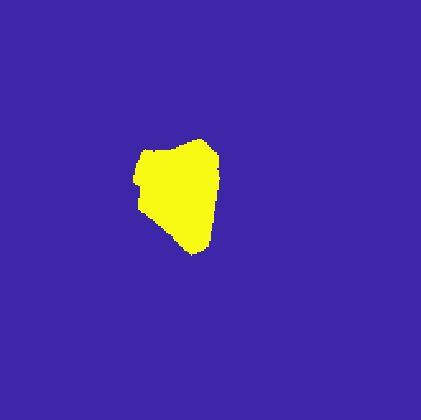}}&\phantom{x}&	\includegraphics[height=3cm]{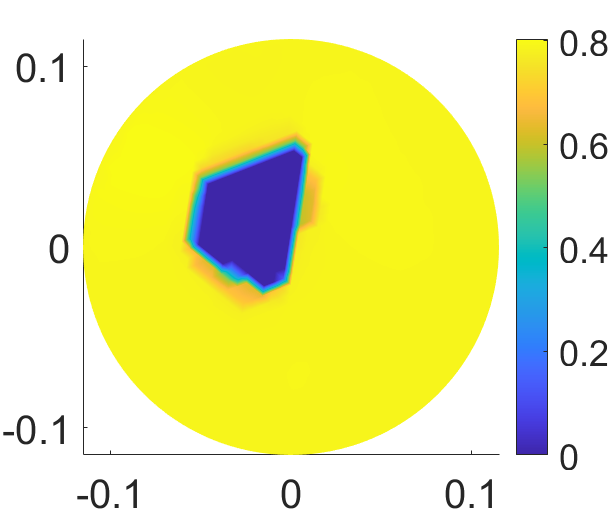}&	\raisebox{0.3cm}{\includegraphics[height=2.5cm]{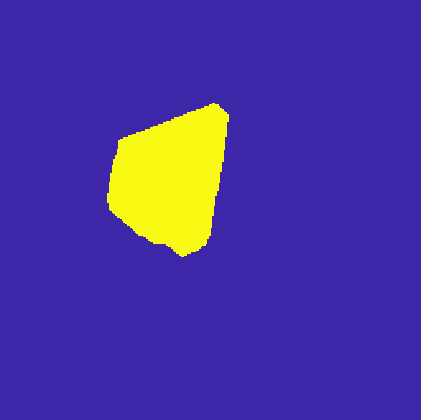}}
	\end{tabular}
\caption{Pairs of output $\sigma$-estimates and corresponding segmentations for phantom \#3 as functions of different numbers $N_{inj}$ of injected currents.}
\label{fig:obj3}
\end{figure}

The performance is consistent even with a smaller number of injections, in particular for $N_{inj}\geq 44$, as shown in Figure \ref{fig:obj3}. For these tests, the hybrid IAS was run with $\eta^{(1)}\in[10^{-5},10^{-4}]$ and $\vartheta^*\in[0.05,0.15]$. The stability in the shape of the reconstructions is also reflected by the SSIM scores that are reported in the fifth column of Table \ref{tab:1}. 

\begin{figure}
	\centering
		\setlength{\tabcolsep}{1pt} 
	\begin{tabular}{ccc}
		phantom \#1&phantom \#2&phantom \#3\\
\includegraphics[height=3.3cm]{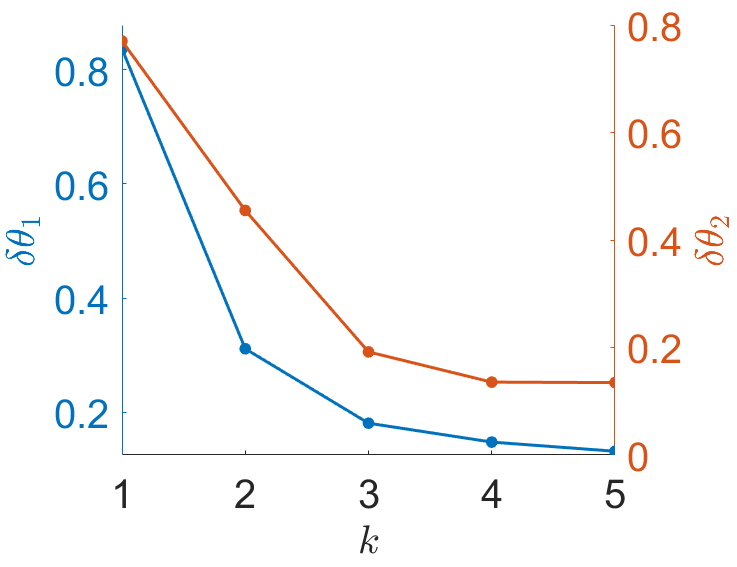}&\includegraphics[height=3.3cm]{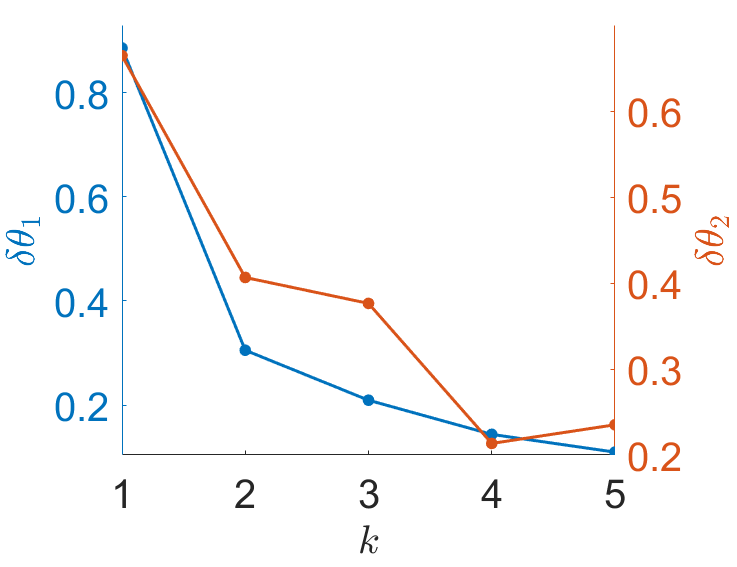}&\includegraphics[height=3.3cm]{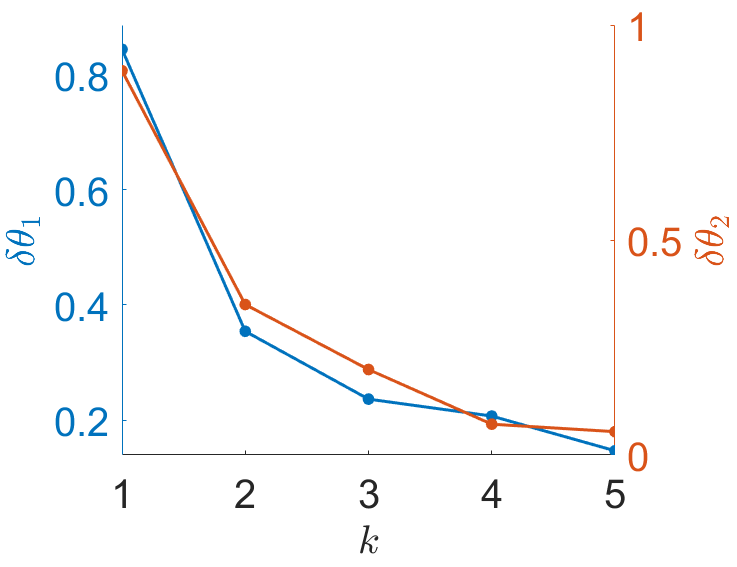}
	\end{tabular}
\caption{For the three phantoms, behavior of  $\delta\theta_1$ and $\delta\theta_2$ as functions of the iteration number of the hybrid IAS in the case of $L=32$ active electrodes and $N_{inj}=76$ current injections.}
\label{fig:rcth}
\end{figure}

As highlighted at the beginning of this section, the stopping criterion selected for the algorithm is only based on the iteration number, as high quality results can be recovered in the very first iterations of the scheme. However, it is also worth analyzing the behavior of the relative changes $\delta\theta_1$, $\delta\theta_2$ defined in \eqref{eq:rcth}, which are shown in Figure \ref{fig:rcth} for the three phantoms for $L=32$ active electrodes, i.e., $N_{inj}=76$ current injections. Notice the overall decreasing behavior, suggesting robust convergence properties of the algorithm.

\begin{figure}
	\centering
	\begin{tabular}{ccc}
				phantom \#1&phantom \#2&phantom \#3\\
	\includegraphics[width=4.0cm]{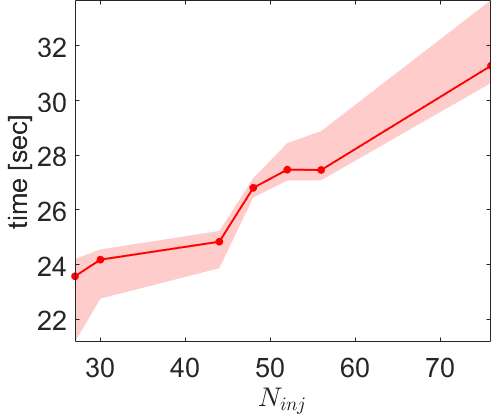}&\includegraphics[width=4.0cm]{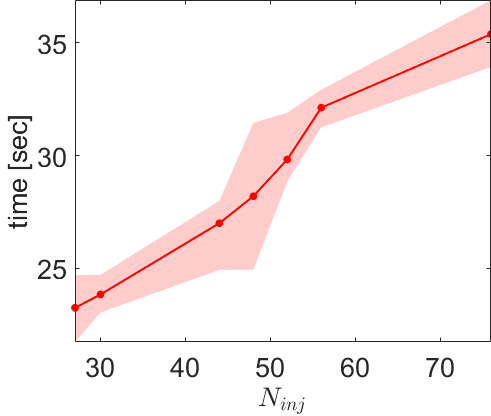}&	\includegraphics[width=4.0cm]{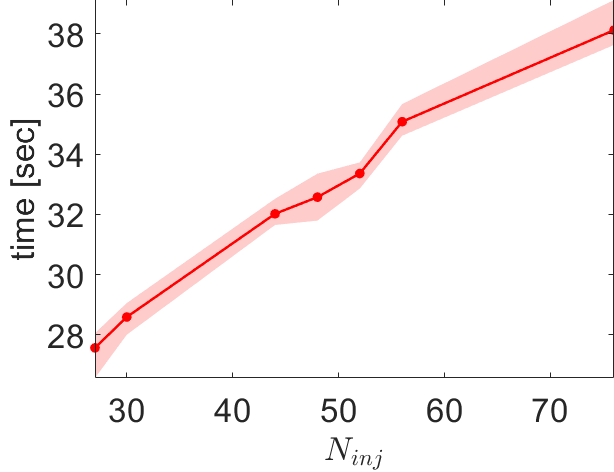}
	\end{tabular}
\caption{For the three phantoms, average running times in seconds and dispersion bands of the hybrid IAS as a function of the number of injected currents.}
\label{fig:time}
\end{figure}

We conclude this section reporting the average computing times over 10 runs of the hybrid IAS for the three phantoms and different numbers of current injections in Figure~\ref{fig:time}, together with the dispersion bands ($\pm$ one standard deviation from the mean). Combining such analysis with the scores reported in Table \ref{tab:1}, one can conclude that high-quality results can be obtained removing up to 8 electrodes from the original setup with a gain in terms of computing times of about 20$\%$. We also highlight that for $L\leq 28$, i.e. $N_{inj}\leq 52$, the linear system solved in steps 3a, 4a, is underdetermined, so that the overall scheme takes advantage of formula \eqref{eq:linsys2}.

\section{Conclusions}\label{sec:concl}

In this work we addressed the solution of the non-linear EIT inverse reconstruction problem with a priori sparsity information. More specifically, after recasting the task in Bayesian terms, we provided a detailed description of the sparsity promoting IAS algorithm when applied in the non-linear settings. A specific attention was given to the parameter selection, and on how the latter can be tailored so as to successfully process the 2023 Kuopio Tomography challenge dataset. Our analysis highlights the efficiency and the flexibility of the discussed algorithm.

\section*{Acknowledgements}

The authors acknowledge the partial support by the NSF, grants DMS 1951446 for Daniela Calvetti and DMS 2204618 for Erkki Somersalo. The support from Solomon R. Guggenheim Foundation to Erkki Somersalo is gratefully acknowledged. Monica Pragliola acknowledges the National Group for Scientific Computation (INdAM-GNCS), Research Projects 2023, and the FRA (Fondi Ricerca Ateneo, University of Naples Federico II) project `HyRED'.

\bibliographystyle{plain}

\bibliography{ref}

\end{document}